\documentclass[10pt, pdftex,a4paper, bibtotoc, 
idxtotoc,
abstracton]{scrartcl}
\usepackage[utf8]{inputenc}
\usepackage[english]{babel}
\usepackage{amssymb,amsmath,amsthm}
\usepackage{natbib}
\setkomafont{sectioning}{\bfseries} 
\usepackage{nccfoots}
\usepackage{csquotes}

\newtheorem{Satz}{Theorem}
\newtheorem*{Satz*}{Theorem A}
\newtheorem{Lem}{Lemma}
\newtheorem{Kor}{Corollary}
\theoremstyle{definition}

\allowdisplaybreaks

\newcommand{\ignore}[1]{}

\makeatletter
\renewcommand*{\@biblabel}[1]{\makebox[\labelwidth][l]{[#1]}}
\makeatother
\DeclareMathOperator{\Var}{Var}

\DeclareMathOperator{\E}{E}

\DeclareMathOperator{\rank}{rank}

\usepackage[svgnames]{xcolor}
\usepackage{listings}
\usepackage[left=25mm,right=25mm,top=25mm,bottom=30mm]{geometry}

\begin{document}

\title{Testing for change-points in long-range dependent time series by means of a self-normalized Wilcoxon test}
\author{Annika Betken\thanks{annika.betken@rub.de}\\ \normalsize{\textit{Fakultät für Mathematik, Ruhr-Universität Bochum, 44780  Bochum,
Germany.}}}
\date{}
\maketitle

\Footnotetext{ }{Research supported by Collaborative Research Center SFB 823 {\em Statistical modeling of nonlinear dynamic processes.}}

\begin{abstract}
We propose a testing procedure based on the Wilcoxon two-sample test
statistic in order to test for  change-points in the mean of long-range dependent data. We show that the corresponding self-normalized test statistic
converges in distribution to a non-degenerate limit under
the hypothesis that no change occurred and that it diverges to infinity under the alternative of a change-point with constant height. Furthermore,  we derive the
asymptotic distribution of the self-normalized Wilcoxon test statistic under local alternatives, that is 
under the assumption that the height of the level shift decreases as the sample size increases.
Regarding the finite sample performance, simulation results confirm that 
 the self-normalized Wilcoxon test yields a consistent discrimination between hypothesis and alternative and that
its
empirical size is  already close to the  significance level 
for moderate sample sizes.\\

\end{abstract}

\section{Introduction}

We consider a data set generated by a stochastic process $(X_i)_{i\geq 1}$,
\begin{align*}
X_i=\mu_i+\varepsilon_i,
\end{align*} 
where $(\mu_i)_{i\geq 1}$ are unknown constants and where $(\varepsilon_i)_{i\geq 1}$ is a stationary, long-range dependent (LRD, in short) process   with mean zero and finite variance.
In particular, we assume that 
\begin{align}
\varepsilon_i=G(\xi_i), \ i\geq 1,
\end{align}
where $(\xi_i)_{i\geq 1}$ is a stationary  Gaussian process
with mean $0$, variance $1$ and long-range dependence, that is  with autocovariance function $\rho$ satisfying
\begin{align*}
\rho(k)\sim k^{-D}L(k), \  k\geq 1, 
\end{align*}
where $0<D< 1$ (referred to as long-range dependence (LRD) parameter)  and where $L$ is a slowly varying function. 
Furthermore, we suppose that $G:\mathbb{R}\longrightarrow \mathbb{R}$ is  a measurable function with $\E\left(G(\xi_i)\right)=0$.

Provided that the previous assumptions hold for the
observations $X_1, \ldots, X_n$,  we wish to test the  hypothesis 
\begin{align*}
H: \mu_1=\ldots =\mu_n
\end{align*}
against the alternative  
\begin{align*}
A: \mu_1=\ldots =\mu_k\neq \mu_{k+1}=\ldots =\mu_n
\end{align*}
for some $k\in \left\{1, \ldots, n-1\right\}$.
Within this setting the location of the change-point is unknown under the alternative.
In order to motivate our choice of a change-point test, we temporarily assume  that 
the change-point location is known, i.e. for a given $k\in \left\{1, \ldots, n-1\right\}$ we consider the alternative
\begin{align*}
A_k: \mu_1=\ldots =\mu_k\neq \mu_{k+1}=\ldots=\mu_n.
\end{align*}
For the test problem $(H, A_k)$,
the Wilcoxon two-sample rank test rejects the hypothesis of no change in the mean for large absolute values of the
test statistic
\begin{align*}
W_{k, n}=\sum\limits_{i=1}^k\sum\limits_{j=k+1}^n
\left(1_{\left\{X_i\leq X_j\right\}}-\frac{1}{2}\right).
\end{align*}
The Wilcoxon change-point test for the test problem $(H, A)$ is defined by reference
to the test statistic $W_{k, n}$; see \cite{DehlingRoochTaqqu2013a}.
It rejects the hypothesis  for large values of
\begin{align*}
\max\limits_{1\leq k\leq n-1}
\left|W_{k, n}\right|=
\max\limits_{1\leq k\leq n-1}
\left|\sum\limits_{i=1}^k\sum\limits_{j=k+1}^n
\left(1_{\left\{X_i\leq X_j\right\}}-\frac{1}{2}\right)\right|.
\end{align*}

With the objective of calculating the asymptotic distribution of the Wilcoxon test statistic under the null hypothesis,  \cite{DehlingRoochTaqqu2013a} consider the stochastic process
\begin{align*}
W_n(\lambda)=\frac{1}{n d_n}\sum\limits_{i=1}^{\lfloor nr\rfloor}\sum\limits_{j=\lfloor nr\rfloor+1}^n\left(1_{\left\{X_i\leq X_j\right\}}-\int_{\mathbb{R}}F(x)dF(x)\right), \ 0\leq \lambda\leq 1,
\end{align*}
where $d_n$ denotes an appropriate normalization.
Assuming that  $(X_i)_{i\geq 1}$ has a continuous marginal distribution function $F$,
the asymptotic distribution of $W_n$
can be derived from the empirical process invariance principle of \cite{DehlingTaqqu1989} as shown in \cite{DehlingRoochTaqqu2013a}. 
It turns out that both, the limit of $W_n$ and the normalization $d_n$, depend on the 
Hermite expansion  
\begin{align*}
1_{\left\{G(\xi_i)\leq x\right\}}-F(x)=\sum\limits_{q=1}^{\infty}\frac{J_q(x)}{q !}H_q(\xi_i),
\end{align*}
where $H_q$ denotes the $q$-th order Hermite polynomial and where
\begin{align*}
J_q(x)=\E \left(H_q(\xi_i)1_{\left\{G(\xi_i)\leq x\right\}}\right).
\end{align*}
The scaling factor $d_n$ is defined by
\begin{align*}
d_n^2=\Var\left(\sum\limits_{j=1}^nH_m(\xi_j)\right),
\end{align*}
where $m$ designates the Hermite rank of the  class of functions $\left\{1_{\left\{G(\xi_i)\leq x\right\}}-F(x), \ x\in \mathbb{R}\right\}$ defined by
\begin{align*}
m:=\min \left\{q\geq 1: J_q(x)\neq 0 \ \text{for  some} \ x\in \mathbb{R}\right\}.
\end{align*}

Presuming the previous conditions hold and  the  long-range dependence parameter $D$ meets the condition $0<D<\frac{1}{m}$,
the process 
\begin{align*}
W_n(\lambda)=\frac{1}{n d_n}\sum\limits_{i=1}^{\lfloor n\lambda\rfloor}\sum\limits_{j=\lfloor n\lambda\rfloor+1}^n\left(1_{\left\{X_i\leq X_j\right\}}-\int_{\mathbb{R}}F(x)dF(x)\right), \ 0\leq \lambda\leq 1,
\end{align*}
converges in distribution to
\begin{align*}
\frac{1}{m!}(Z_m(\lambda)-\lambda Z_m(1))\int_{\mathbb{R}}J_m(x)dF(x), \ 0\leq \lambda \leq 1,
\end{align*}
where $\left(Z_m(\lambda)\right)_{\lambda \in \left[0, 1\right]}$ is an $m$-th order Hermite process, which is self-similar with parameter $H=1-\frac{mD}{2}\in \left(\frac{1}{2}, 1\right)$.
If $m=1$, the Hermite process $Z_m$ 
 equals a standard fractional Brownian motion process with Hurst parameter $H=1-\frac{D}{2}$. We refer 
to \cite{Taqqu1979} for a general definition of the Hermite process $Z_m$.

An application of the continuous mapping theorem to the process $W_n$ yields the asymptotic distribution of the Wilcoxon change-point test. 
More precisely, it has been proved by \cite{DehlingRoochTaqqu2013a} that under the hypothesis of no change in the mean,
the Wilcoxon test statistic
\begin{align*}
\frac{1}{n d_n}\max\limits_{1\leq k\leq n-1}\left|\sum\limits_{i=1}^k\sum\limits_{j=k+1}^n\left(1_{\left\{X_i\leq X_j\right\}}-\frac{1}{2}\right)\right|
\end{align*}
converges in distribution to
\begin{align*}
\sup\limits_{0\leq \lambda \leq 1}\left|\frac{1}{m!}(Z_m(\lambda)-\lambda Z_m(1))\right|\left|\int_{\mathbb{R}}J_m(x)dF(x)\right|.
\end{align*}

Furthermore,
\cite{DehlingRoochTaqqu2013b}
 investigate the asymptotic behaviour of the Wilcoxon change-point test under  the alternative
with the objective of determining  the height of the level shift in such a way that the power of the self-normalized Wilcoxon test is non-trivial.  
For this purpose, they consider
 local alternatives defined by
\begin{align*}
A_{\tau, h_n}:
\mu_i=\begin{cases}
\mu \ &\text{for} \ i=1, \ldots, \lfloor n\tau\rfloor\\
\mu+h_n \ &\text{for} \ i=\lfloor n\tau \rfloor +1, \ldots, n,
\end{cases}
\end{align*}
where $0< \tau < 1$ and where $h_n\sim c\frac{d_n}{n}$, so that under the sequence of local alternatives $A_{\tau, h_n}$ the height of the level shift decreases if the sample size increases.
Under
the additional assumption that 
$G(\xi_i)$
has a continuous distribution function
$F$
with bounded
density
$f$, this guarantees that under the sequence of alternatives $A_{\tau, h_n}$, the process
\begin{align*}
\frac{1}{nd_n}\sum\limits_{i=1}^{\lfloor n\lambda\rfloor}\sum\limits_{j=\lfloor n\lambda\rfloor +1}^n\left(1_{\{X_i\leq X_j\}}-\frac{1}{2}\right), \ 0\leq \lambda\leq 1,
\end{align*}
converges in distribution to the limit process
\begin{align*}
\frac{1}{m!}(Z_m(\lambda)-\lambda Z_m(1))\int_{\mathbb{R}}J_m(x)dF(x)+c\delta_{\tau}(\lambda)\int_{\mathbb{R}}f^2(x)dx, \  0\leq \lambda\leq 1,
\end{align*}
where $\delta_{\tau}:[0, 1]\longrightarrow \mathbb{R}$  is defined by
\begin{align*}
\delta_{\tau}(\lambda)=
\begin{cases}
\lambda(1-\tau) \ &\text{for} \ \lambda\leq \tau\\
(1-\lambda)\tau \ &\text{for} \ \lambda\geq \tau
\end{cases}.
\end{align*}
By another application of the continuous mapping theorem 
it then follows that  the Wilcoxon change-point test converges in distribution to a non-degenerate limit process under the sequence of local alternatives $A_{\tau, h_n}$; see \cite{DehlingRoochTaqqu2013b}.

\vspace{5mm}

\section{Main Results}

An application of
 the Wilcoxon change-point test to a given data set
 presupposes that the scaling factor $d_n$
is known. 
Usually this is not the case in statistical practice so that in general the Wilcoxon change-point test as proposed in \cite{DehlingRoochTaqqu2013a} depends on an unknown normalization.
As an alternative we propose a normalization that only depends on the given  realizations and therefore is referred to as self-normalization.
The self-normalization approach we consider
 has originally been established  in another context; see \cite{Lobato2001}.
It has been
 extended   to the change-point testing problem by  \cite{ShaoZhang2010} in order to test for change-points in the mean of short-range dependent time series. These authors used the self-normalization method on the Kolmogorov-Smirnov test statistic, in doing so also taking the change-point alternative into account.
Lobato  as well as Shao and Zhang considered weak dependent processes only. Following the approach  in Shao and Zhang an application to possibly long-range dependent processes was introduced by Shao, who established a self-normalized version of the CUSUM change-point test; see  \cite{Shao2011}.

As the CUSUM test has the disadvantage of not being robust against possible outliers in the data, an extension of the self-normalization idea to the Wilcoxon test statistic leads to a change-point test that not only has the advantage of avoiding the choice of unknown parameters but also yields a robust alternative 
to the CUSUM test.

Given observations $X_1, \ldots, X_n$, we consider the rank statistics  defined by
\begin{flalign*}
R_i=\rank(X_i)=\sum\limits_{j=1}^n1_{\{X_j\leq X_i\}} 
\end{flalign*}
for $i=1, \ldots, n$.
An extension of the self-normalization approach to the Wilcoxon change-point test is based on an application of the CUSUM change-point test  in terms of the rank statistics $R_i$.
Note that due to the identity
\begin{align*}
\max\limits_{k}\left|\sum\limits_{i=1}^kR_i-\frac{k}{n}\sum\limits_{i=1}^nR_i\right|
=\max\limits_{k}\left|\sum\limits_{i=1}^k\sum\limits_{j=k+1}^n\left(1_{\left\{X_i\leq X_j\right\}}-\frac{1}{2}\right)\right|,
\end{align*}
 the CUSUM test statistic of the ranks equals the Wilcoxon change-point test statistic.
Instead of dividing the test statistic (which is the maximum taken among every possible outcome of the Wilcoxon two-sample rank test) by the unknown quantity $nd_n$ we consider a normalization  factor
that depends on the location of a potential change-point and which therefore is different for every possible outcome of the Wilcoxon two-sample rank test.

 We define
\begin{align*}
G_n(k)
=\frac{\left|\sum\limits_{i=1}^kR_i-\frac{k}{n}\sum\limits_{i=1}^nR_i\right|}{
\left\{\frac{1}{n}\sum\limits_{t=1}^k S_t^2(1,k)+\frac{1}{n}\sum\limits_{t=k+1}^n S_t^2(k+1,n)\right\}^{\frac{1}{2}}
}, 
\end{align*}
where 
\begin{align*}
&S_{t}(j, k)=\sum\limits_{h=j}^t\left(R_h-\bar{R}_{j, k}\right),\\
&\bar{R}_{j, k}=\frac{1}{k-j+1}\sum\limits_{t=j}^kR_t.
\end{align*}

The self-normalized Wilcoxon test  rejects the hypothesis $H:\mu_1=\ldots =\mu_n$ for
large values of the test statistic 
\begin{align*}
T_n(\tau_1, \tau_2)=\sup_{k\in \left[\lfloor n\tau_1\rfloor, \lfloor n\tau_2\rfloor\right]}G_n(k), 
\end{align*}
where $0< \tau_1 <\tau_2 <1$.

Note that the proportion of the data that is included in the calculation of the supremum 
 is restricted by the choice of $\tau_1$ and $\tau_2$.
 This is important as 
the choice of $\tau_1$ and $\tau_2$ influences the properties of the test.
Structural breaks  at the
beginning or the end of a sample are hard to detect since there is a lack of information concerning the behaviour of the time series before or after a 
potential break point.
 Hence, the interval $\left[\tau_1, \tau_2\right]$ must be small enough for the critical values not to get too large on the one hand, yet large enough to include potential break points on the other hand. A common  choice is $\tau_1= 1-\tau_2=0.15$; see \cite{Andrews1993}.

The following theorem
states the asymptotic  distribution of the test statistic $T_n(\tau_1, \tau_2)$ 
under the hypothesis of no change in the mean.

\begin{Satz}\label{asymptotic  distribution under H}
Suppose that $\left(X_i\right)_{i\geq 1}$ is a stationary process with continuous distribution function $F$  defined by
\begin{align*}
X_i=\mu_i+G(\xi_i)
\end{align*}
for unknown constants $\left(\mu_i\right)_{i\geq 1}$ and a stationary, long-range dependent Gaussian  process  $\left(\xi_i\right)_{i\geq 1}$  with mean $0$, variance $1$ and LRD parameter $0<D <\frac{1}{m}$, where 
 $m$ denotes the Hermite rank of the class of functions $1_{\left\{G(\xi_i)\leq x\right\}}-F(x)$, $x \in \mathbb{R}$.  
Moreover, assume that   $\int_{\mathbb{R}}J_m(x)dF(x)\neq 0$ and that $G:\mathbb{R}\longrightarrow \mathbb{R}$ is a measurable function.
Then, under the hypothesis of no
change in the mean, it follows that $T_n(\tau_1, \tau_2)~\overset{\mathcal{D}}{\longrightarrow}~T(m,\tau_1, \tau_2)$, where
\begin{align*}
&T(m, \tau_1, \tau_2)=\sup\limits_{\lambda\in \left[\tau_1, \tau_2\right]}\frac{ \left|Z_m(\lambda)-\lambda Z_m(1)\right|}{
\Bigl\{\int_0^{\lambda}\left(V_m(r; 0, \lambda)\right)^2dr+\int_{\lambda}^1 \left(V_m(r; \lambda, 1)\right)^2dr\Bigr\}^{\frac{1}{2}}}
\intertext{with}
&V_m(r; r_1, r_2)=Z_m(r)-Z_m(r_1)-\frac{r-r_1}{r_2-r_1}\left\{Z_m(r_2)-Z_m(r_1)\right\}
\end{align*}
for $r\in \left[r_1, r_2\right]$, $0< r_1< r_2< 1$.
\end{Satz}

As  consistency under fixed alternatives is considered as a fundamental characteristic of appropriate hypothesis testing, we aim at proving  Theorem \ref{consistency}, which implies that if there is a change-point in the mean of constant height, the empirical power of the self-normalized Wilcoxon test tends to $1$.
For this purpose, we suppose that under the alternative
\begin{align}\label{change-point model}
X_i=\begin{cases}
\mu +G(\xi_i), \ i=1, \ldots, k^*, \\
\mu +\Delta +G(\xi_i),  \ i=k^*+1, \ldots, n,
\end{cases}
\end{align}
where $k^*=\lfloor n\tau\rfloor$ and  $\Delta \neq 0$
is fixed.

\begin{Satz}\label{consistency}
Suppose that $(\xi_i)_{i\geq 1}$  is a stationary, long-range dependent Gaussian process with mean $0$, variance $1$ and LRD parameter $D$. Moreover, let  $G:\mathbb{R}\longrightarrow \mathbb{R}$ be a measurable function and assume  that $G(\xi_i)$ has a continuous distribution function $F$.
Given that the parameter $D$  satisfies $0< D<\frac{1}{m}$, where
 $m$ denotes the Hermite rank of the class of functions  $1_{\{G(\xi_i)\leq x\}}-F(x)$, $x\in \mathbb{R}$,
$T_n(\tau_1, \tau_2)$ diverges in probability to $\infty$  under fixed alternatives, i.e. if $\left(X_i\right)_{i\geq 1}$ satisfies \eqref{change-point model}.
\end{Satz}

Furthermore, we wish to study the asymptotic behaviour  of the self-normalized Wilcoxon change-point test under local alternatives defined by
\begin{align*}
A_{\tau, h_n}(n):
\mu_i=\begin{cases}
\mu \ &\text{for} \ i=1, \ldots, \lfloor n\tau\rfloor,\\
\mu+h_n \ &\text{for} \ i=\lfloor n\tau \rfloor +1, \ldots, n,
\end{cases}
\end{align*}
where $0< \tau < 1$ and $h_n\longrightarrow 0$. 
The following theorem confirms that the self-normalized Wilcoxon test statistic converges to a non-degenerate limit under the sequence of local alternatives  $A_{\tau, h_n}$.
\begin{Satz}\label{asymptotic distribution under local A}
Suppose that $(\xi_i)_{i\geq 1}$  is a stationary Gaussian process with mean $0$, variance $1$ and autocovariance function  
\begin{align*}
\rho(k)\sim k^{-D}L(k),
\end{align*}
where $L$  is a slowly varying function and where $0< D<\frac{1}{m}$. Moreover, let  $G:\mathbb{R}\longrightarrow \mathbb{R}$ be a measurable function. We assume  that $G(\xi_i)$ has a continuous distribution function $F$ with bounded density $f$. Let $m$ denote the Hermite rank of the class of functions  $1_{\{G(\xi_i)\leq x\}}-F(x)$, $x\in \mathbb{R}$, and suppose that   $\int_{\mathbb{R}}J_m(x)dF(x)\neq 0$. Then, under the sequence of alternatives  $A_{\tau, h_n}$ with  $h_n\sim c\frac{d_n}{n}$, it follows that $T_n(\tau_1, \tau_2)$
converges in distribution to
\begin{align*}
T(m, \tau_1, \tau_2)
=\sup\limits_{\lambda\in \left[\tau_1, \tau_2\right]}\frac{\left|\frac{1}{m!}\int_{\mathbb{R}}J_m(x)dF(x)(Z_m(\lambda)-\lambda Z_m(1))+c\delta_{\tau}(\lambda)\int_{\mathbb{R}}f^2(x)dx\right|}{\left\{\int_0^\lambda\left(V_{m, \tau}(r; 0, \lambda)\right)^2dr+\int_{\lambda}^1\left(V_{m, \tau}(r; \lambda, 1)\right)^2dr\right\}^{\frac{1}{2}}},
\end{align*}
where
\begin{align*}
V_{m, \tau}(r; 0, \lambda)
&=\frac{1}{m!}\int_{\mathbb{R}}J_m(x)dF(x)\left(Z_m(r)-\frac{r}{\lambda}Z_m(\lambda)\right)+c\int_{\mathbb{R}}f^2(x)dx\left(\delta_{\tau}(r)-\frac{r}{\lambda}\delta_{\tau}(\lambda)\right),\\
V_{m, \tau}(r; \lambda, 1)&=\frac{1}{m!}\int_{\mathbb{R}}J_m(x)dF(x)\left\{Z_m(r)-Z_m(\lambda)-\frac{r-\lambda}{1-\lambda}\left(Z_m(1)-Z_m(\lambda)\right)\right\}\\
&\quad \ +c\int_{\mathbb{R}}f^2(x)dx
\left(\delta_{\tau}(r)-\frac{1-r}{1-\lambda}\delta_{\tau}(\lambda)\right).
\end{align*}
\end{Satz}

\vspace{5mm}

\section{Simulation studies}

We will now investigate the finite sample performance of the self-normalized Wilcoxon test statistic.
For this purpose, we take 
$G(t) = t$
so that $\left(X_i\right)_{i\geq 1}$ is a Gaussian process. 
Since $G$ is strictly increasing, the Hermite coefficient $J_1(x)$ is not equal to $0$ for all $x\in \mathbb{R}$; see \cite{DehlingRoochTaqqu2013a}.   Therefore, it holds that $m = 1$, where $m$ denotes the Hermite rank of $1_{\left\{G(\xi_i)\leq x\right\}}-F(x), x\in \mathbb{R}$.  As a result, $T_n(\tau_1, \tau_2)$
has approximately the same distribution as
\begin{align*}
&\sup\limits_{\lambda\in \left[\tau_1, \tau_2\right]}\frac{ \left|B_H(\lambda)-\lambda B_H(1)\right|}{
\Bigl\{\int_0^{\lambda}\left(V_H(r; 0, \lambda)\right)^2dr+\int_{\lambda}^1 \left(V_H(r; \lambda, 1)\right)^2dr\Bigr\}^{\frac{1}{2}}}\\
&V_H(r; r_1, r_2)=B_H(r)-B_H(r_1)-\frac{r-r_1}{r_2-r_1}\left\{B_H(r_2)-B_H(r_1)\right\},
\end{align*}
where $B_H$ is a fractional Brownian motion process with Hurst parameter $H=1-\frac{D}{2}$.

We set critical values on the basis of
 $10,000$
simulations of fractional Brownian motion time series for different Hurst parameters $H$ and different levels of significance; see Table \ref{critical values SN-Wilcoxon}. 

\begin{table}[h]
\center			
\begin{tabular}{r r r r}
 & 10\%	& 5\% 	& 1\% 	\\	
\hline			
$H$ = 0.6 & 6.182835 	& 7.276568	& 9.785915\\	
$H$ = 0.7 & 6.847260 	& 8.190125	& 11.380584\\
$H$ = 0.8 & 7.767277	& 9.495194	& 13.021080\\
$H$ = 0.9 & 8.520039	& 10.333602	& 14.544094	
\end{tabular}
\caption{Simulated critical values for the distribution of $T(1, \tau_1, \tau_2)$ when  $\left[\tau_1, \tau_2\right]=\left[0.15, 0.85\right]$. The sample size is  $1 000$, the number of
replications is $10, 000$.} 
\label{critical values SN-Wilcoxon}	
\end{table}

The calculation of the relative frequency of false rejections under the hypothesis is based on 
 $10,000$ realizations of fractional Gaussian noise time series with varying length; see Table \ref{level SN-Wilcoxon}.

\begin{table}[htbp]
\center	
\begin{tabular}{r c c c c c }
n & & \multicolumn{1}{c}{H=0.6} &  \multicolumn{1}{c}{H=0.7} &  \multicolumn{1}{c}{H=0.8} &  \multicolumn{1}{c}{H=0.9} \\
\hline
10 & & 0.057 & 0.052  & 0.036  & 0.026 \\
50 & & 0.048 & 0.050 & 0.046 & 0.052 
\\
100 & & 0.049 & 0.055 & 0.050 & 0.053 \\
500 & &  0.053 & 0.050 &  0.049 & 0.054 \\
1000 & & 0.053 & 0.053 &0.050 & 0.052 
\end{tabular}
\caption{Level of the self-normalized Wilcoxon change-point test for fractional Gaussian noise time series of length $n$ with Hurst parameter $H$.
The level of significance is $5\%$. The calculations are based on $10,000$ simulation
runs.}
\label{level SN-Wilcoxon}	
\end{table}

The simulation results suggest that the self-normalized Wilcoxon test performs well under the hypothesis since empirical size and asymptotic significance level are already close for moderate sample sizes.
In particular, it is  notable
that the size of the self-normalized Wilcoxon change-point test differs considerably from the size of the original Wilcoxon change-point test 
when $H=0.9$, that means when we
have very strong dependence. In that case, the convergence of the Wilcoxon change-point test statistic 
appears to be rather slow under the hypothesis (see \cite{DehlingRoochTaqqu2013a}, Table 2), whereas the size  of the self-normalized Wilcoxon change-point test is 
still close to the corresponding  level of significance.

We consider fractional Gaussian noise time series with a level shift of height $\Delta$
after a proportion $\tau$ of the data  in order  to analyse the behaviour of the test statistic under the alternative. We have done so for several choices of  $\Delta$ and $\tau$ and for sample sizes $n=100$ and $n = 500$.

\begin{table}
\center		
\begin{tabular}{l l l l l l l l l l l l }
& & & \multicolumn{2}{c}{$\Delta$ = 0.5} & & \multicolumn{2}{c}{$\Delta$ = 1}
& & \multicolumn{2}{c}{$\Delta$ = 2} \\
\cline{4-5}\cline{7-8}\cline{10-11}
& & & 10\% 	& 5\% & 	 
	& 10\% 	& 5\% & 	 
	& 10\% 	& 5\% &	
	\\	
\hline	
$H$ = 0.6 & $n$ = 100& 
&0.474	&0.348&	&0.956	&0.916&	&1.000	&1.000
\\
 		   & $n$ = 500& 
&0.941	&0.898&	&1.000	&1.000&	&1.000	&1.000 		   	
\\							
$H$ = 0.7 & $n$ = 100& 
&0.355	&0.239&	&0.801	&0.696& &1.000	&0.999
\\
 & $n$ = 500& 
&0.655	&0.530&	&0.992	&0.985&	&1.000	&1.000
\\					 		
$H$ = 0.8 & $n$ = 100& 
&0.281	&0.177&	&0.652	&0.526&	&0.993	&0.984
\\
 & $n$ = 500& 
&0.405	&0.297&	&0.872	&0.802&	&1.000	&1.000
\\							
$H$ = 0.9 & $n$ = 100& 
&0.426	&0.287&	&0.744	&0.628&	&0.992	&0.983
\\
 & $n$ = 500& 
&0.412	&0.280&	&0.786	&0.689&	&0.998	&0.997
\end{tabular}
\caption{Empirical power of the self-normalized Wilcoxon change-point test
for fractional Gaussian noise of length $n=100$ and $n=500$ with Hurst parameter $H $ and a level shift in the mean of height $\Delta$ after a proportion $\tau=0.5$. The calculations are based on $5,000$ simulation runs.} 
\label{Power SN-Wilcoxon tau=0.5}
\end{table}

\begin{table}[htbp]
\center			
\begin{tabular}{l l l l l l l l l l l l }
& & & \multicolumn{2}{c}{$\Delta$ = 0.5} & & \multicolumn{2}{c}{$\Delta$ = 1}
& & \multicolumn{2}{c}{$\Delta$ = 2} \\
\cline{4-5}\cline{7-8}\cline{10-11}

& & & 10\% 	& 5\% & 	 
	& 10\% 	& 5\% & 	 
	& 10\% 	& 5\% &	
	\\	
\hline	
$H$ = 0.6 & $n$ = 100& 
&0.321	&0.204&	&0.813	&0.690&	&1.000	&1.000
\\
 & $n$ = 500& 
&0.795	&0.678&	&1.000	&0.999&	&1.000	&1.000
\\							
$H$ = 0.7 & $n$ = 100& 
&0.222	&0.125&	&0.570	&0.401&	&0.989	&0.968
\\
 & $n$ = 500& 
&0.437	&0.309&	&0.948	&0.891&	&1.000	&1.000
\\					 		
$H$ = 0.8 & $n$ = 100& 
&0.195	&0.106&	&0.417	&0.265 &	&0.931	&0.839
\\
 & $n$ = 500& 
&0.264	&0.164&	&0.682	&0.530&	&0.999	&0.995
\\							
$H$ = 0.9 & $n$ = 100& 
&0.339	&0.198&	&0.578	&0.403&	&0.961	&0.889
\\
 & $n$ = 500& 
&0.312	&0.186&	&0.612	&0.442&	&0.989	&0.966
\end{tabular}
\caption{Empirical power of the self-normalized Wilcoxon change-point test
for fractional Gaussian noise of length $n=100$ and $n=500$ with Hurst parameter $H $ and a level shift in the mean of height $\Delta$ after a proportion $\tau=0.25$. The calculations are based on $5,000$ simulation runs.} 
\label{Power SN-Wilcoxon tau=0.25}
\end{table}

The simulations of the empirical power confirm that the rejection rate becomes higher when  $\Delta$ increases. 
 Comparing the empirical power for different Hurst parameters $H$, we note that the test tends to have less power as $H$ becomes large. 
This seems natural since when there is very strong dependence, i.e. $H$
is large, the variance of the series  increases, so that it becomes  harder to detect a level shift of a fixed height.
 In addition, change-points that are located in the middle of the sample are detected more often than change-points that are located close   to the boundary of the testing region determined by $\left[\tau_1, \tau_2\right]$.
Furthermore, Table \ref{Power SN-Wilcoxon tau=0.25} and Table \ref{Power SN-Wilcoxon tau=0.5} show that  an increasing sample size goes along with an increase of the empirical power.
This result confirms that the self-normalized Wilcoxon change-point test yields a consistent discrimination between hypothesis and alternative.

\newpage

\section{Proofs}

In order to simplify notation, we write 
\begin{align*}
&J(x)=\frac{1}{m!}J_m(x),\\
&Z(\lambda)=Z_m(\lambda).
\end{align*}

\textit{Proof of Theorem \ref{asymptotic  distribution under H}.} The essential step in the proof
of Theorem \ref{asymptotic  distribution under H} is to find a representation for the test  statistic $T_n(\tau_1, \tau_2)$ as a functional of the Wilcoxon process
\begin{align*}
W_n(\lambda)=\frac{1}{nd_n}\sum\limits_{i=1}^{\lfloor n\lambda\rfloor}\sum\limits_{j=\lfloor n\lambda\rfloor+1}^n\left(1_{\left\{X_i\leq X_j\right\}}-\frac{1}{2}\right), \ 0\leq \lambda\leq 1.
\end{align*}
For this purpose,  rewrite
\begin{align*}
G_n(k)&=\frac{\left|\sum\limits_{i=1}^{k}\sum\limits_ {j=k+1}^{n}\left(1_{\{X_i\leq X_j\}}-\frac{1}{2}\right)\right|}{
\left\{\frac{1}{n}\sum\limits_{t=1}^k S_t^2(1,k)+\frac{1}{n}\sum\limits_{t=k+1}^n S_t^2(k+1,n)\right\}^{\frac{1}{2}}
}\\
&=\frac{\frac{1}{nd_n}\left|\sum\limits_{i=1}^{k}\sum\limits_ {j=k+1}^{n}\left(1_{\{X_i\leq X_j\}}-\frac{1}{2}\right)\right|}{
\frac{1}{nd_n}\left\{\frac{1}{n}\sum\limits_{t=1}^k S_t^2(1,k)+\frac{1}{n}\sum\limits_{t=k+1}^n S_t^2(k+1,n)\right\}^{\frac{1}{2}}
}.
\end{align*}
As we have
\begin{align*}
\frac{1}{nd_n}\left|\sum\limits_{i=1}^{k}\sum\limits_ {j=k+1}^{n}\left(1_{\{X_i\leq X_j\}}-\frac{1}{2}\right)\right|=\left|W_n(\lambda)\right|
\end{align*}
for the numerator of $G_n(k)$ if $k=\lfloor n\lambda\rfloor$,
it remains to show that
the denominator of $G_n(k)$ can be represented as a functional of $W_n$.
Since
\begin{align*}
R_i=n+1-\sum\limits_{j=1}^n1_{\left\{X_i\leq X_j\right\}} 
\end{align*}
almost surely,
it follows that
\begin{align*}
S_t(1, k)
=&-\sum\limits_{h=1}^t\left(\sum\limits_{j=1}^n1_{\{X_h\leq X_j\}}-\frac{1}{k}\sum\limits_{i=1}^k\sum\limits_{j=1}^n1_{\{X_i\leq X_j\}}\right)\\
=&-\Biggl\{\sum\limits_{i=1}^t\sum\limits_{j=t+1}^n\left(1_{\{X_i\leq X_j\}}-\frac{1}{2}\right)+\sum\limits_{i=1}^t\sum\limits_{j=1}^t\left(1_{\{X_i\leq X_j\}}-\frac{1}{2}\right)\\
&-\frac{t}{k}\sum\limits_{i=1}^k\sum\limits_{j=k+1}^n\left(1_{\{X_i\leq X_j\}}-\frac{1}{2}\right)
-\frac{t}{k}\sum\limits_{i=1}^k\sum\limits_{j=1}^k\left(1_{\{X_i\leq X_j\}}-\frac{1}{2}\right)\Biggr\}
\end{align*}
almost surely.
Moreover, it is well known that  
\begin{align}\label{indicator}
\sum\limits_{i=1}^l\sum\limits_{j=1}^l1_{\left\{X_i\leq X_j\right\}}=\frac{l(l+1)}{2}.
\end{align}
Hence,
\begin{align*}
\sum\limits_{i=1}^l\sum\limits_{j=1}^l\left(1_{\left\{X_i\leq X_j\right\}}-\frac{1}{2}\right)=\frac{l(l+1)}{2}-\frac{l^2}{2}=\frac{l}{2},
\end{align*}
so that
\begin{align*}
S_t(1, k)
&=-\Biggl\{\sum\limits_{i=1}^t\sum\limits_{j=t+1}^n\left(1_{\{X_i\leq X_j\}}-\frac{1}{2}\right)+\frac{t}{2}
-\frac{t}{k}\sum\limits_{i=1}^k\sum\limits_{j=k+1}^n\left(1_{\{X_i\leq X_j\}}-\frac{1}{2}\right)
-\frac{t}{k}\frac{k}{2}\Biggr\}\\
&=-\Biggl\{\sum\limits_{i=1}^t\sum\limits_{j=t+1}^n\left(1_{\{X_i\leq X_j\}}-\frac{1}{2}\right)
-\frac{t}{k}\sum\limits_{i=1}^k\sum\limits_{j=k+1}^n\left(1_{\{X_i\leq X_j\}}-\frac{1}{2}\right)\Biggr\}
\end{align*}
almost surely.
Thus, if  $\lambda\in \left[\tau_1, \tau_2\right]$,
\begin{align*}
&\int_0^{\lambda} \left(\sum\limits_{i=1}^{\lfloor nr\rfloor}\sum\limits_{j=\lfloor nr\rfloor+1}^n\left(1_{\{X_i\leq X_j\}}-\frac{1}{2}\right)
-\frac{\lfloor nr\rfloor}{\lfloor n\lambda\rfloor}\sum\limits_{i=1}^{\lfloor n\lambda\rfloor}\sum\limits_{j=\lfloor n\lambda\rfloor+1}^n\left(1_{\{X_i\leq X_j\}}-\frac{1}{2}\right)\right)^2dr\\
&=\sum\limits_{t=0}^{\lfloor n\lambda\rfloor}\int_{\frac{t}{n}}^{\frac{t+1}{n}}\left(\sum\limits_{i=1}^{\lfloor nr\rfloor}\sum\limits_{j=\lfloor nr\rfloor+1}^n\left(1_{\{X_i\leq X_j\}}-\frac{1}{2}\right)
-\frac{\lfloor nr\rfloor}{\lfloor n\lambda\rfloor}\sum\limits_{i=1}^{\lfloor n\lambda\rfloor}\sum\limits_{j=\lfloor n\lambda\rfloor+1}^n\left(1_{\{X_i\leq X_j\}}-\frac{1}{2}\right)\right)^2dr\\
&\quad \ -\int_{\lambda}^{\frac{\lfloor n\lambda\rfloor+1}{n}}\left(\sum\limits_{i=1}^{\lfloor nr\rfloor}\sum\limits_{j=\lfloor nr\rfloor+1}^n\left(1_{\{X_i\leq X_j\}}-\frac{1}{2}\right)
-\frac{\lfloor nr\rfloor}{\lfloor n\lambda\rfloor}\sum\limits_{i=1}^{\lfloor n\lambda\rfloor}\sum\limits_{j=\lfloor n\lambda\rfloor+1}^n\left(1_{\{X_i\leq X_j\}}-\frac{1}{2}\right)\right)^2dr,
\end{align*}
where
\begin{align*}
\sum\limits_{i=1}^{\lfloor nr\rfloor}\sum\limits_{j=\lfloor nr\rfloor+1}^n\left(1_{\{X_i\leq X_j\}}-\frac{1}{2}\right)
-\frac{\lfloor nr\rfloor}{\lfloor n\lambda\rfloor}\sum\limits_{i=1}^{\lfloor n\lambda\rfloor}\sum\limits_{j=\lfloor n\lambda\rfloor+1}^n\left(1_{\{X_i\leq X_j\}}-\frac{1}{2}\right)=0
\end{align*}
for $r\in \left[\lambda, \frac{\lfloor n\lambda\rfloor+1}{n}\right)$.
Therefore, the integral over  that interval
equals $0$.
Consequently, 
\begin{align*}
&\int_0^{\lambda} \left(\sum\limits_{i=1}^{\lfloor nr\rfloor}\sum\limits_{j=\lfloor nr\rfloor+1}^n\left(1_{\{X_i\leq X_j\}}-\frac{1}{2}\right)
-\frac{\lfloor nr\rfloor}{\lfloor n\lambda\rfloor}\sum\limits_{i=1}^{\lfloor n\lambda\rfloor}\sum\limits_{j=\lfloor n\lambda\rfloor+1}^n\left(1_{\{X_i\leq X_j\}}-\frac{1}{2}\right)\right)^2dr\\
&=\sum\limits_{t=0}^{\lfloor n\lambda\rfloor}\int_{\frac{t}{n}}^{\frac{t+1}{n}}\left(\sum\limits_{i=1}^{\lfloor nr\rfloor}\sum\limits_{j=\lfloor nr\rfloor+1}^n\left(1_{\{X_i\leq X_j\}}-\frac{1}{2}\right)
-\frac{\lfloor nr\rfloor}{\lfloor n\lambda\rfloor}\sum\limits_{i=1}^{\lfloor n\lambda\rfloor}\sum\limits_{j=\lfloor n\lambda\rfloor+1}^n\left(1_{\{X_i\leq X_j\}}-\frac{1}{2}\right)\right)^2dr\\
&=\frac{1}{n}\sum\limits_{t=0}^{k}\left(\sum\limits_{i=1}^{t}\sum\limits_{j= t+1}^n\left(1_{\{X_i\leq X_j\}}-\frac{1}{2}\right)
-\frac{t}{k}\sum\limits_{i=1}^{k}\sum\limits_{j=k+1}^n\left(1_{\{X_i\leq X_j\}}-\frac{1}{2}\right)\right)^2\\
&=\frac{1}{n}\sum\limits_{t=1}^{k} S_t^2(1,k)
\end{align*}
almost surely in case $k=\lfloor n\lambda\rfloor$.

For the second term in the denominator of $G_n(k)$
the following equations hold almost surely
\begin{align*}
S_t(k+1, n)
=&-\Biggr\{\sum\limits_{h=k+1}^t\left(\sum\limits_{j=1}^n1_{\{X_h\leq X_j\}}-\frac{1}{n-k}\sum\limits_{i=k+1}^n\sum\limits_{j=1}^n1_{\{X_i\leq X_j\}}\right)\Biggl\}\\
=&-\Biggr\{\sum\limits_{i=1}^t\sum\limits_{j=t+1}^n\left(1_{\{X_i\leq X_j\}}-\frac{1}{2}\right)+\sum\limits_{i=1}^t\sum\limits_{j=1}^t\left(1_{\{X_i\leq X_j\}}-\frac{1}{2}\right)\\
&-\sum\limits_{i=1}^k\sum\limits_{j=k+1}^n\left(1_{\{X_i\leq X_j\}}-\frac{1}{2}\right)-
\sum\limits_{i=1}^k\sum\limits_{j=1}^k\left(1_{\{X_i\leq X_j\}}-\frac{1}{2}\right)\\
&-\frac{t-k}{n-k}\sum\limits_{i=k+1}^n\sum\limits_{j=1}^k\left(1_{\{X_i\leq X_j\}}-\frac{1}{2}\right)-
\frac{t-k}{n-k}\sum\limits_{i=k+1}^n\sum\limits_{j=k+1}^n\left(1_{\{X_i\leq X_j\}}-\frac{1}{2}\right)\Biggl\}.
\end{align*}
By \eqref{indicator} we get
\begin{align*}
\sum\limits_{i=k+1}^n\sum\limits_{j=k+1}^n\left(1_{\{X_i\leq X_j\}}-\frac{1}{2}\right)=\frac{(n-k)(n-k+1)}{2}-\frac{(n-k)^2}{2}=\frac{n-k}{2}.
\end{align*}
Furthermore,
\begin{align*}
1_{\{X_i\leq X_j\}}-\frac{1}{2}=1-1_{\{X_j< X_i\}}-\frac{1}{2}=-\left(1_{\{X_j\leq X_i\}}-\frac{1}{2}\right)
\end{align*}
almost surely if $i\neq j$.   
This yields
\begin{align*}
S_t(k+1, n)
=&-\Biggr\{\sum\limits_{i=1}^t\sum\limits_{j=t+1}^n\left(1_{\{X_i\leq X_j\}}-\frac{1}{2}\right)+\frac{t}{2}
-\sum\limits_{i=1}^k\sum\limits_{j=k+1}^n\left(1_{\{X_i\leq X_j\}}-\frac{1}{2}\right)-
\frac{k}{2}\\
&-\frac{t-k}{n-k}\sum\limits_{i=k+1}^n\sum\limits_{j=1}^k\left(1_{\{X_i\leq X_j\}}-\frac{1}{2}\right)-
\frac{t-k}{n-k}\frac{n-k}{2}\Biggl\}\\
=&-\Biggr\{\sum\limits_{i=1}^t\sum\limits_{j=t+1}^n\left(1_{\{X_i\leq X_j\}}-\frac{1}{2}\right)
-\frac{n-t}{n-k}\sum\limits_{i=1}^k\sum\limits_{j=k+1}^n\left(1_{\{X_i\leq X_j\}}-\frac{1}{2}\right)\Biggl\}.
\end{align*}
We obtain  for $\lambda\in \left[\tau_1, \tau_2\right]$
\begin{align*}
&\int_{\lambda}^1 \left(\sum\limits_{i=1}^{\lfloor nr\rfloor}\sum\limits_{j=\lfloor nr\rfloor+1}^n\left(1_{\{X_i\leq X_j\}}-\frac{1}{2}\right)
-\frac{n-\lfloor nr\rfloor}{n-\lfloor n\lambda\rfloor}\sum\limits_{i=1}^{\lfloor n\lambda\rfloor}\sum\limits_{j=\lfloor n\lambda\rfloor+1}^n\left(1_{\{X_i\leq X_j\}}-\frac{1}{2}\right)\right)^2dr\\
&=\sum\limits_{t=\lfloor n\lambda\rfloor+1}^{n-1}\int_{\frac{t}{n}}^{\frac{t+1}{n}}\left(\sum\limits_{i=1}^{\lfloor nr\rfloor}\sum\limits_{j=\lfloor nr\rfloor+1}^n\left(1_{\{X_i\leq X_j\}}-\frac{1}{2}\right)
-\frac{n-\lfloor nr\rfloor}{n-\lfloor n\lambda\rfloor}\sum\limits_{i=1}^{\lfloor n\lambda\rfloor}\sum\limits_{j=\lfloor n\lambda\rfloor+1}^n\left(1_{\{X_i\leq X_j\}}-\frac{1}{2}\right)\right)^2dr\\
&\quad \ +\int_{\lambda}^{\frac{\lfloor n\lambda\rfloor+1}{n}}\left(\sum\limits_{i=1}^{\lfloor nr\rfloor}\sum\limits_{j=\lfloor nr\rfloor+1}^n\left(1_{\{X_i\leq X_j\}}-\frac{1}{2}\right)
-\frac{n-\lfloor nr\rfloor}{n-\lfloor n\lambda\rfloor}\sum\limits_{i=1}^{\lfloor n\lambda\rfloor}\sum\limits_{j=\lfloor n\lambda\rfloor+1}^n\left(1_{\{X_i\leq X_j\}}-\frac{1}{2}\right)\right)^2dr
\end{align*}
almost surely, where
\begin{align*}
\sum\limits_{i=1}^{\lfloor nr\rfloor}\sum\limits_{j=\lfloor nr\rfloor+1}^n\left(1_{\{X_i\leq X_j\}}-\frac{1}{2}\right)
-\frac{n-\lfloor nr\rfloor}{n-\lfloor n\lambda\rfloor}\sum\limits_{i=1}^{\lfloor n\lambda\rfloor}\sum\limits_{j=\lfloor n\lambda\rfloor+1}^n\left(1_{\{X_i\leq X_j\}}-\frac{1}{2}\right)
=0
\end{align*}
if $r\in \left[\lambda, \frac{\lfloor n\lambda\rfloor+1}{n}\right)$. Therefore, the integral over that interval equals $0$. 
For $k=\lfloor n\lambda\rfloor$ this implies
\begin{align*}
&\int_{\lambda}^1 \left(\sum\limits_{i=1}^{\lfloor nr\rfloor}\sum\limits_{j=\lfloor nr\rfloor+1}^n\left(1_{\{X_i\leq X_j\}}-\frac{1}{2}\right)
-\frac{n-\lfloor nr\rfloor}{n-\lfloor n\lambda\rfloor}\sum\limits_{i=1}^{\lfloor n\lambda\rfloor}\sum\limits_{j=\lfloor n\lambda\rfloor+1}^n\left(1_{\{X_i\leq X_j\}}-\frac{1}{2}\right)\right)^2dr\\
&=\frac{1}{n}\sum\limits_{t=k +1}^{n-1}\left(\sum\limits_{i=1}^{t}\sum\limits_{j=t+1}^n\left(1_{\{X_i\leq X_j\}}-\frac{1}{2}\right)
-\frac{n-t}{n-k}\sum\limits_{i=1}^{k}\sum\limits_{j=k+1}^n\left(1_{\{X_i\leq X_j\}}-\frac{1}{2}\right)\right)^2\\
&=\frac{1}{n}\sum\limits_{t=k+1}^{n-1}S_{t}^2(k+1,n)\\
&=\frac{1}{n}\sum\limits_{t=k+1}^{n}S_{t}^2(k+1,n).
\end{align*}

Due to the previous considerations, the properly normalized denominator of $G_n(k)$ can (almost surely) be represented as follows
\begin{align*}
&\frac{1}{nd_n}\left\{\frac{1}{n}\sum\limits_{t=1}^{k}S_t^2(1, k)+\frac{1}{n}\sum\limits_{t=k+1}^n S_t^2(k+1,n)\right\}^\frac{1}{2}\\
&=\left\{\int_0^{\lambda}\left(W_n(r)-\frac{c_n(r)}{c_n(\lambda)}W_n(\lambda)\right)^2dr+\int_{\lambda}^1\left(W_n(r)-\frac{1-c_n(r)}{1-c_n(\lambda)} W_n(\lambda)\right)^2dr\right\}^{\frac{1}{2}}, 
\end{align*}
where $c_n(\lambda)=\frac{\lfloor n\lambda\rfloor}{n}$ for $\lambda\in \left[0, 1\right]$.
All in  all, this yields
\begin{align*}
T_n(\tau_1, \tau_2)
=\sup\limits_{\lambda\in \left[\tau_1, \tau_2\right]}\frac{\left|W_n(\lambda)\right|}{\left\{\int_0^{\lambda}\left(W_n(r)-\frac{c_n(r)}{c_n(\lambda)}W_n(\lambda)\right)^2dr+\int_{\lambda}^1\left(W_n(r)-\frac{1-c_n(r)}{1-c_n(\lambda)} W_n(\lambda)\right)^2dr\right\}^{\frac{1}{2}}}.
\end{align*}

The foregoing characterization of the self-normalized Wilcoxon test statistic  points out that 
a representation  of $T_n(\tau_1, \tau_2)$ as a functional of the process
\begin{align*}
W_n(\lambda)=\frac{1}{n d_n}\sum\limits_{i=1}^{\lfloor n\lambda\rfloor}\sum\limits_{j=\lfloor n\lambda\rfloor+1}^n\left(1_{\left\{X_i\leq X_j\right\}}-\frac{1}{2}\right), \ 0\leq \lambda\leq 1,
\end{align*}
also depends on the function series $\left(c_n\right)_{n\in \mathbb{N}}$ in $D\left[0, 1\right]$ defined by
$c_n(\lambda)=\frac{\lfloor n\lambda\rfloor}{n}, \ 0\leq \lambda \leq 1$.
Since
\begin{align*}
\sup\limits_{\lambda\in \left[0, 1\right]}\left|\frac{\lfloor n\lambda\rfloor}{n}-\lambda\right|
=\sup\limits_{\lambda\in \left[0, 1\right]}\left(\lambda-\frac{\lfloor n\lambda\rfloor}{n}\right)
\leq \sup\limits_{\lambda\in \left[0, 1\right]}\left(\lambda-\frac{ n\lambda-1}{n}\right)
=\frac{1}{n}\longrightarrow 0,
\end{align*}
the sequence $c_n$, $n\in \mathbb{N}$, converges  with respect to the supremum norm  to $c\in D\left[0, 1\right]$ defined by 
$c(\lambda)=\lambda$ for $\lambda\in \left[0, 1\right]$.
To simplify subsequent calculations, we treat $c_n$ and $c$ as random variables with values in the closure of
\begin{align*}
M=\left\{f\in D\left[0, 1\right]\left|\right. f(\lambda)=\frac{\lfloor n\lambda\rfloor}{n} \ \text{for some } n\in \mathbb{N}, \ n\geq \frac{1}{\tau_1}\right\}.
\end{align*} 
Note that 
\begin{align*}
h_n=\left(\begin{array}{c}
c_n\\
W_n
\end{array}\right)
\overset{\mathcal{D}}{\longrightarrow}
\left(\begin{array}{c}
c\\
W_m^*
\end{array}\right),
\end{align*}
where
\begin{align}
W_m^*(\lambda)=(Z(\lambda)-\lambda Z(1))\int_{\mathbb{R}}J(x)dF(x), \ 0\leq \lambda \leq 1.
\end{align}

Obviously,
the self-normalized Wilcoxon test statistic can be represented as a functional of the random vector $h_n$.
Hence, an application of the continuous  mapping theorem just requires the definition
of an appropriate function $G:\overline{M}\times D\left[0, 1\right]\longrightarrow \mathbb{R}$ that maps $h_n$ on $T_n(\tau_1, \tau_2)=G(h_n)$. 
For $\lambda\in \left[\tau_1, \tau_2\right]$ consider the function $G_{\lambda}:\overline{M}\times D\left[0, 1\right]\longrightarrow \mathbb{R}$ that maps an element $h=\left(h_1, h_2\right)$ on
\begin{align*}
 \frac{\left|
 h_2(\lambda)\right|}{\left\{\int_0^{\lambda}\left(h_2(r)-\frac{h_1(r)}{h_1(\lambda)}h_2(\lambda)\right)^2dr+\int_{\lambda}^1\left(h_2(r)-\frac{1 -h_1(r)}{1-h_1(\lambda)}h_2(\lambda)\right)^2dr\right\}^{\frac{1}{2}}},
\end{align*}
provided 
 that the function $F:\overline{M}\times D\left[0, 1\right]\longrightarrow \mathbb{R}$ defined by
\begin{align*}
F(h)=  
\inf\limits_{\lambda\in \left[\tau_1, \tau_2\right]}\left\{\int_0^{\lambda}\left(h_2(r)-\frac{h_1(r)}{h_1(\lambda)}h_2(\lambda)\right)^2dr
+\int_{\lambda}^1\left(h_2(r)-\frac{1 -h_1(r)}{1-h_1(\lambda)}h_2(\lambda)\right)^2dr\right\}^{\frac{1}{2}}
\end{align*}
does not equal $0$ in $h$.
Given that  $h\in F^{-1}\left(\left\{0\right\}\right)$, we set $G_{\lambda}(h)=-1$.

Since $T_n(\tau_1, \tau_2)=\sup_{\lambda\in \left[\tau_1, \tau_2\right]}G_{\lambda}(h_n)$, 
we intend to apply the continuous mapping theorem to the function
$G: \overline{M}\times D\left[0, 1\right]\longrightarrow\mathbb{R}$, where $G(h)=\sup_{\lambda\in \left[\tau_1, \tau_2\right]}G_{\lambda}(h)$.
Thus, we have to 
verify that the  function $G$ complies with the requirements of the continuous mapping theorem, i.e. we have to prove the following assertions:
\begin{enumerate}
 \item[1)] The function $G$ is measurable with respect to the uniform product metric on $\overline{M}\times D\left[0, 1\right]$.
\item[2)] We have $P(h \in D_G)=0$, where $D_G$ denotes the set of discontinuities of $G$.
\end{enumerate}

In order to show that $G$ is measurable, we consider the restrictions of $G$ to   $\left(\overline{M}\times D\left[0, 1\right]\right)\setminus F^{-1}(\{0\})$ and $F^{-1}(\{0\})$,
respectively.
Both restrictions are continuous with respect to the uniform metric. In particular, both restrictions are Borel measurable. Since the restricted domains are Borel measurable subsets of $\overline{M}\times D\left[0, 1\right]$, the measurability of the restrictions implies the measurability of $G$.

It remains to show that $P(h \in D_G)=0$.
Again,  consider the restriction of $G$ to $\left(\overline{M}\times D\left[0, 1\right]\right)\setminus F^{-1}(\{0\})$.
Because of the continuity of the restriction, 
 $G$ is continuous at every $h\in \left(\overline{M}\times D\left[0, 1\right]\right)\setminus F^{-1}(\{0\})$ as $F^{-1}(\{0\})$
is a closed subset of $\overline{M}\times D\left[0, 1\right]$.
Therefore, $D_G$ is a subset of $F^{-1}(\{0\})$.   Consequently, it suffices to show that $P(h\in F^{-1}(\{0\}))=0$ in order to prove that $P(h\in D_G)=0$.

The random vector $h=(c, W_m^*)$ is an element of $F^{-1}(\{0\})$ if and only if 
the expression
\begin{align}\label{inf}
\inf\limits_{\lambda\in \left[\tau_1, \tau_2\right]}\left\{\int_0^{\lambda}\left(f(r)-\frac{r}{\lambda}f(\lambda)\right)^2dr+\int_{\lambda}^1\left(f(r)-\frac{1-r}{1-\lambda}f(\lambda)\right)^2dr\right\}^{\frac{1}{2}}
\end{align}
vanishes when $f=W_m^*$.

Note that
\begin{align}\label{equation1}
W_m^*(r)-\frac{r}{\lambda}W_m^*(\lambda)
&=\int J(x)dF(x)\left\{\left(Z(r)-rZ(1)\right)-\frac{r}{\lambda}\left(Z(\lambda)-\lambda Z(1)\right)\right\} \notag\\
&=\int J(x)dF(x)\left\{Z(r)-\frac{r}{\lambda}Z(\lambda)\right\}
\end{align}
and
\begin{align}\label{equation2}
W_m^*(r)-\frac{1-r}{1-\lambda}W_m^*(\lambda)
&=\int J(x)dF(x)\left\{\left(Z(r)-rZ(1)\right)-\frac{1-r}{1-\lambda}\left(Z(\lambda)-\lambda Z(1)\right)\right\} \notag\\
&=\int J(x)dF(x)\left\{Z(r)-Z(\lambda)-\frac{r-\lambda}{1-\lambda}\left(Z(1)-Z(\lambda)\right)\right\}.
\end{align}
Therefore, and as $Z \in C\left[0, 1\right]$ almost surely (see \cite{MaejimaTudor2007}),  the term in formula \eqref{inf} vanishes if for some  $\lambda\in \left[\tau_1, \tau_2\right]$
\begin{align*}
\left\{\int_0^{\lambda}\left(V_m(r;  0, \lambda)\right)^2dr+\int_{\lambda}^1\left(V_m(r; \lambda, 1)\right)^2dr\right\}^{\frac{1}{2}}=0,
\end{align*}
where
\begin{align*}
V_m(r; r_1, r_2)=Z(r
)-Z(r_1)-\frac{r-r_1}{r_2-r_1}\left(Z(r_2)-Z(r_1)\right).
\end{align*}
It suffices to show that the sample paths of $W_m^*$ do not belong to  the set of continuous functions $f$ that satisfy 
\begin{align}\label{equation3}
\left\{\int_0^{\lambda}\left(f(r)-\frac{r}{\lambda}f(\lambda)\right)^2dr+\int_{\lambda}^1\left(f(r)-f(\lambda)-\frac{r-\lambda}{1-\lambda}\left(f(1)-f(\lambda)\right)\right)^2dr\right\}^{\frac{1}{2}}=0
\end{align}
for some $\lambda\in \left[\tau_1, \tau_2\right]$.
The above equation only holds if the integrands vanish almost surely on the corresponding intervals.
In particular, a continuous function $f \in D\left[0, 1\right]$  that meets formula \eqref{equation3} 
satisfies
\begin{align*}
f(r)=\frac{1}{\lambda}f(\lambda)r
\end{align*}
if $r\in\left[0, \lambda\right]$
and
\begin{align*}
f(r)&=f(\lambda)+\frac{r-\lambda}{1-\lambda}\left\{f(1)-f(\lambda)\right\}\\
&=f(\lambda)-\frac{\lambda}{1-\lambda}\left\{f(1)-f(\lambda)\right\}+\frac{1}{1-\lambda}\left\{f(1)-f(\lambda)\right\}r
\end{align*}
if $r\in \left[\lambda, 1\right]$. Consequently, the set of continuous functions which lie in $F^{-1}(\{0\})$ corresponds to the class of functions
\begin{align*}
A = & \Bigl\{f\in D\left[0, 1\right]\left|\right. \  \text{for some $\lambda\in [\tau_1, \tau_2]$ and $a, b \in \mathbb{R}$}\\
&f(r)=\frac{1}{\lambda}a r \ \text{on $[0, \lambda]$} \ \text{ and} \\
&f(r)=a-\frac{\lambda}{1-\lambda}\left\{b-a\right\}+\frac{1}{1-\lambda}\left\{b-a\right\}r \ \text{on $[\lambda, 1]$}\Bigr\}.
\end{align*}
It follows that $P(Z\in A)=0$ because
the sample paths of the Hermite process $Z$ are nowhere differentiable with probability $1$ (see \cite{Mikosch1998}), whereas an element in $A$ is differentiable almost everywhere.
This implies $P(h\in D_G)=0$. 
 
Having verified the preconditions of the continuous mapping theorem
we are now able to conclude that the test statistic
 $T_n(\tau_1, \tau_2)$ converges in distribution to
\begin{align*}
T(m,\tau_1, \tau_2)=\sup\limits_{\lambda\in \left[\tau_1, \tau_2\right]}\frac{\left|W_m^*(\lambda)\right|}{\left\{\int_0^{\lambda}\left(W_m^*(r)-\frac{r}{\lambda}W_m^*(\lambda)\right)^2dr+\int_{\lambda}^1\left(W_m^*(r)-\frac{1-r}{1-\lambda} W_m^*(\lambda)\right)^2dr\right\}^{\frac{1}{2}}}.
\end{align*}

Due to  \eqref{equation1} and  \eqref{equation2}, the limit process $T(m, \tau_1, \tau_2)$ equals  
\begin{align*}
&\sup\limits_{\lambda\in \left[\tau_1, \tau_2\right]}\frac{ \left|Z(\lambda)-\lambda Z(1)\right|}{
\Bigl\{\int_0^{\lambda}\left(V_m(r; 0, \lambda)\right)^2dr+\int_{\lambda}^1 \left(V_m(r; \lambda, 1)\right)^2dr\Bigr\}^{\frac{1}{2}}}.
\end{align*}
Thus, we have established Theorem 
\ref{asymptotic  distribution under H}.
\hfill $ \Box$

\vspace{5mm}

In the proof of  Theorem \ref{consistency} we make use of  preliminary results stated in Lemma \ref{Lemma1}, Lemma \ref{Lemma2} and Corollary \ref{Corollary1}.
The line of argument  that  verifies Lemma \ref{Lemma1} and Lemma \ref{Lemma2}
is a modification of the proof that establishes Theorem 3.1  in \cite{DehlingRoochTaqqu2013b}.

\begin{Lem}\label{Lemma1}
Suppose that $\left(\xi_i\right)_{i\geq 1}$  is a stationary, long-range dependent Gaussian  process   with mean $0$, variance $1$ and LRD parameter $0<D <\frac{1}{m}$, where 
 $m$ denotes the Hermite rank of the class of functions $1_{\left\{G(\xi_i)\leq x\right\}}-F(x)$, $x \in \mathbb{R}$.  
Moreover, assume that $\left(
G(\xi_i)\right)_{i\geq 1}$ has a  continuous distribution function $F$ and that $G:\mathbb{R}\longrightarrow \mathbb{R}$ is a measurable function. Then, if $\Delta \in \mathbb{R}$,
\begin{align*}
&\frac{1}{n^2}\sum\limits_{i=1}^{\lfloor n\lambda\rfloor}\sum\limits_{j=\lfloor n\tau\rfloor +1}^n1_{\left\{G(\xi_i)\leq G(\xi_j)+\Delta\right\}}\overset{P}{\longrightarrow}
\lambda(1-\tau)\int_{\mathbb{R}}F(x+\Delta)dF(x),\\
&\frac{1}{n^2}\sum\limits_{i=1}^{\lfloor n\tau\rfloor}\sum\limits_{j=\lfloor n\lambda\rfloor +1}^n1_{\left\{G(\xi_i)\leq G(\xi_j)+\Delta\right\}}
\overset{P}{\longrightarrow}
\tau(1-\lambda)\int_{\mathbb{R}}F(x+\Delta)dF(x)
\end{align*}
for fixed $\tau$, uniformly in $\lambda\leq \tau$ and $\lambda\geq \tau$, respectively.
\end{Lem}

\textit{Proof of Lemma \ref{Lemma1}.}
We give a proof for the first assertion only as
the convergence of the second term follows by an analogous argumentation.

Let $F_k$ and $F_{k+1, n}$ denote the empirical distribution  functions of the first $k$ and last $n-k$ realizations of $G(\xi_1), \ldots, G(\xi_n)$,  i.e.
\begin{align*}
&F_k(x)=\frac{1}{k}\sum\limits_{i=1}^k1_{\left\{G(\xi_i)\leq x\right\}},\\
&F_{k+1, n}(x)=\frac{1}{n-k}\sum\limits_{i=k+1}^n1_{\left\{G(\xi_i)\leq x\right\}}.
\end{align*}
For $\lambda\leq \tau$ this yields the following representation: 
\begin{align*}
\sum\limits_{i=1}^{\lfloor n\lambda\rfloor}\sum\limits_{j=\lfloor n\tau\rfloor+1}^n1_{\left\{G(\xi_i)\leq G(\xi_j)+\Delta\right\}}
&=\left(n-\lfloor n\tau\rfloor\right)\lfloor n\lambda\rfloor
\frac{1}{n-\lfloor n\tau\rfloor}\sum\limits_{j=\lfloor n\tau\rfloor +1}^nF_{\lfloor n\lambda\rfloor}( G(\xi_j)+\Delta)
\\
&=\left(n-\lfloor n\tau\rfloor\right)\lfloor n\lambda\rfloor
\int_{\mathbb{R}}F_{\lfloor n\lambda\rfloor}( x+\Delta)
dF_{\lfloor n \tau \rfloor +1, n}(x)
\end{align*}

Since $\frac{n-\lfloor n\tau\rfloor}{n}\longrightarrow 1-\tau$, it suffices to show that $\lfloor n\lambda\rfloor
\int_{\mathbb{R}}F_{\lfloor n\lambda\rfloor}( x+\Delta)
dF_{\lfloor n \tau \rfloor +1, n}(x)$ converges to $\lambda\int_{\mathbb{R}}F(x+\Delta)dF(x)$.
For this purpose, we consider the inequality 
\begin{align}\label{inequality1 Lemma1}
&\sup\limits_{0\leq \lambda\leq \tau}\left|\frac{1}{n}\lfloor n\lambda\rfloor\int_{\mathbb{R}}F_{\lfloor n\lambda\rfloor}(x+\Delta)dF_{\lfloor n\tau\rfloor +1, n}(x)-\lambda\int_{\mathbb{R}}F(x+\Delta)dF(x)\right|
\\
&\leq 
\sup\limits_{0\leq \lambda\leq \tau}
\left|\frac{1}{n}\int_{\mathbb{R}}\lfloor n\lambda\rfloor F_{\lfloor n\lambda\rfloor}(x+\Delta)dF_{\lfloor n\tau\rfloor +1, n}(x)-\frac{\lfloor n\lambda\rfloor}{n}\int_{\mathbb{R}}F(x+\Delta)dF_{\lfloor n\tau\rfloor +1, n}(x)\right|\notag\\
&\quad \ +\sup\limits_{0\leq \lambda\leq \tau}\left|\frac{\lfloor n\lambda\rfloor}{n}\int_{\mathbb{R}}F(x+\Delta)dF_{\lfloor n\tau\rfloor +1, n}(x)
-\lambda \int_{\mathbb{R}}F(x+\Delta)dF_{\lfloor n\tau\rfloor +1, n}(x)
\right|\notag\\
&\quad \ +\sup\limits_{0\leq \lambda\leq \tau}\left|\lambda \int_{\mathbb{R}}F(x+\Delta)dF_{\lfloor n\tau\rfloor +1, n}(x)
-\lambda \int_{\mathbb{R}}F(x+\Delta)dF(x)
\right| \notag
\end{align}
and we will show 
that each of the three terms on its right-hand side converges to $0$.

For the third summand 
we get
\begin{align*}
&\sup\limits_{0\leq \lambda\leq \tau}\left|\lambda \int_{\mathbb{R}}F(x+\Delta)dF_{\lfloor n\tau\rfloor +1, n}(x)
-\lambda \int_{\mathbb{R}}F(x+\Delta)dF(x)
\right|\\
&=\sup\limits_{0\leq \lambda\leq \tau}\left|\lambda\left(1- \int_{\mathbb{R}}F_{\lfloor n\tau\rfloor +1, n}(x)dF(x+\Delta)
- \left(1-\int_{\mathbb{R}}F(x)dF(x+\Delta)\right)
\right)\right|\\
&=\tau\left| \int_{\mathbb{R}}\left(F_{\lfloor n\tau\rfloor +1, n}(x)- F(x)\right)dF(x+\Delta)
\right|\\
&\leq\tau\sup\limits_{x\in \mathbb{R}}\left| F_{\lfloor n\tau\rfloor +1, n}(x)- F(x)\right|
\end{align*}
as a consequence of integration by parts. Furthermore, we have
\begin{align*}
\sup\limits_{x\in \mathbb{R}}\left|F_n(x)-F(x)\right|\longrightarrow 0  \ a.s.
\end{align*}
by an application of the Glivenko-Cantelli theorem (see  \cite{Krengel1985})
 to the stationary and ergodic process $\left(G(\xi_i)\right)_{i\geq 1}$. 
So as to deduce an 	analogous result for $F_{\lfloor n\tau\rfloor +1, n}$ we rewrite
\begin{align*}
F_{\lfloor n\tau\rfloor +1, n}(x)=\frac{n}{n-\lfloor n\tau\rfloor}F_n(x)- \frac{\lfloor n\tau\rfloor}{n-\lfloor n\tau\rfloor}F_{\lfloor n\tau\rfloor}(x)
\end{align*}
and we may therefore conclude
\begin{align*}
\left|F_{\lfloor n\tau\rfloor+1, n}(x)-F(x)\right|
&\leq\left|\frac{n}{n-\lfloor n\tau\rfloor}\right|\left|F_{n}(x)-F(x)\right|
+\left|\frac{\lfloor n\tau\rfloor}{n-\lfloor n\tau\rfloor}\right|\left|F_{\lfloor n\tau\rfloor}(x)-F(x)\right|.
\end{align*}
Thus, 
\begin{align}\label{Glivenko-Cantelli}
\sup\limits_{x\in \mathbb{R}}\left|F_{\lfloor n\tau\rfloor +1, n}(x)-F(x)\right|\longrightarrow 0 \ a.s.
\end{align}
which implies that the third term on the right-hand side of \eqref{inequality1 Lemma1} converges to $0$ almost surely.

Regarding the second term on the right-hand side of \eqref{inequality1 Lemma1}, we obtain
\begin{align*}
&\sup\limits_{0\leq \lambda\leq \tau}\left|\frac{\lfloor n\lambda\rfloor}{n}\int_{\mathbb{R}}F(x+\Delta)dF_{\lfloor n\tau\rfloor +1, n}(x)
-\lambda \int_{\mathbb{R}}F(x+\Delta)dF_{\lfloor n\tau\rfloor +1, n}(x)
\right|\\
&=\sup\limits_{0\leq \lambda\leq \tau}\left|\frac{\lfloor n\lambda\rfloor}{n}-\lambda\right|\left|\int_{\mathbb{R}}F(x+\Delta)dF_{\lfloor n\tau\rfloor +1, n}(x)
\right|.
\end{align*}
The right-hand side of this equation converges to $0$ since $\left|\int_{\mathbb{R}}F(x+\Delta)dF_{\lfloor n\tau\rfloor +1, n}(x)
\right|$ is bounded by $1$, and as 
\begin{align*}
\sup\limits_{0\leq \lambda \leq \tau}\left|\frac{\lfloor n\lambda\rfloor}{n}-\lambda\right|
\longrightarrow 0.
\end{align*}
In order to show that the first term in \eqref{inequality1 Lemma1} converges to $0$ as well, we consider the following inequality:
\begin{align}\label{inequality2 Lemma1}
&\sup\limits_{0\leq \lambda\leq \tau}
\left|\frac{1}{n}\int_{\mathbb{R}}\lfloor n\lambda\rfloor F_{\lfloor n\lambda\rfloor}(x+\Delta)dF_{\lfloor n\tau\rfloor +1, n}(x)-\frac{\lfloor n\lambda\rfloor}{n}\int_{\mathbb{R}}F(x+\Delta)dF_{\lfloor n\tau\rfloor +1, n}(x)\right|\\
&=
\sup\limits_{0\leq \lambda\leq \tau}
\left|\frac{1}{n}\int_{\mathbb{R}}\lfloor n\lambda\rfloor \left(F_{\lfloor n\lambda\rfloor}(x+\Delta)-F(x+\Delta)\right)dF_{\lfloor n\tau\rfloor +1, n}(x)\right|\notag\\
&\leq
\sup\limits_{0\leq \lambda\leq \tau}
\left|\frac{d_n}{n}\int_{\mathbb{R}}d_n^{-1}\lfloor n\lambda\rfloor\left(F_{\lfloor n\lambda\rfloor}(x+\Delta)-F(x+\Delta)\right)-J(x+\Delta)Z(\lambda)dF_{\lfloor n\tau\rfloor +1, n}(x)\right|\notag\\
&\quad \ +\sup\limits_{0\leq \lambda\leq \tau}
\left|\frac{d_n}{n}Z(\lambda)\int_{\mathbb{R}}J(x+\Delta)dF_{\lfloor n\tau\rfloor +1, n}(x)\right|\notag
\end{align}

In what follows, we will prove that both terms
 on the right-hand side of \eqref{inequality2 Lemma1}  converge to $0$.
For this purpose, we make use of the empirical process 
non-central limit theorem
of
\cite{DehlingTaqqu1989} which states that
\begin{align*}
\left(d_n^{-1}\lfloor n\lambda\rfloor (F_{\lfloor n\lambda\rfloor}(x)-F(x))\right)_{x\in \left[-\infty, \infty\right], \lambda \in \left[0, 1\right]}
\overset{\mathcal{D}}{\longrightarrow}J(x)Z(\lambda),
\end{align*}
where     \enquote{$\overset{\mathcal{D}}{\longrightarrow}$} denotes convergence in distribution with respect to the $\sigma$-field generated by the open balls in $D\left(\left[-\infty, \infty\right]\times \left[0, 1\right]\right)$, equipped with the supremum norm.

Due to the Dudley-Wichura version of Skorohod's representation theorem (see \cite{ShorackWellner1986}, Theorem 2.3.4), we may assume 
without loss of generality that
\begin{align*}
\sup\limits_{\lambda, x}\left|d_n^{-1}\lfloor n\lambda\rfloor\left(F_{\lfloor n\lambda\rfloor}(x)-F(x)\right)-J(x)Z(\lambda)\right|\longrightarrow 0
\end{align*}
almost surely; see \cite{DehlingRoochTaqqu2013a}.
As a consequence,
the first summand in \eqref{inequality2 Lemma1} converges to $0$ since
\begin{align*}
&\sup\limits_{0\leq \lambda\leq \tau}
\left|\frac{d_n}{n}\int_{\mathbb{R}}d_n^{-1}\lfloor n\lambda\rfloor\left(F_{\lfloor n\lambda\rfloor}(x+\Delta)-F(x+\Delta)\right)-J(x+\Delta)Z(\lambda)dF_{\lfloor n\tau\rfloor +1, n}(x)\right|\\
&=\frac{d_n}{n}\sup\limits_{0\leq \lambda\leq \tau, x}\left|d_n^{-1}\lfloor n\lambda\rfloor\left(F_{\lfloor n\lambda\rfloor}(x+\Delta)-F(x+\Delta)\right)-J(x+\Delta)Z(\lambda)\right|
\end{align*}
and as $\frac{d_n}{n}$ converges to $0$ as well.

For the second summand we get the following inequality:
\begin{align*}
\sup\limits_{0\leq \lambda\leq \tau}
\left|\frac{d_n}{n}Z(\lambda)\int_{\mathbb{R}}J(x+\Delta)dF_{\lfloor n\tau\rfloor +1, n}(x)\right|
\leq \frac{d_n}{n}\sup\limits_{0\leq \lambda \leq \tau}\left|Z(\lambda)\right|
\left|\int_{\mathbb{R}}J(x+\Delta)dF_{\lfloor n\tau\rfloor +1, n}(x)\right|
\end{align*}
Note that
\begin{align*}
J(x)&=\int_{\mathbb{R}}1_{\left\{G(y)\leq x\right\}}H_m(y)\varphi(y)dy\\
&=\int_{\mathbb{R}}H_m(y)\varphi(y)dy-\int_{\mathbb{R}}1_{\left\{x\leq G(y)\right\}}H_m(y)\varphi(y)dy\\
&=-\int_{\mathbb{R}}1_{\left\{x\leq G(y)\right\}}H_m(y)\varphi(y)dy,
\end{align*}
where $\varphi$  denotes the standard normal density function,
since
\begin{align*}
\int_{\mathbb{R}}H_m(y)\varphi(y)dy=0.
\end{align*}
For this reason, we have
\begin{align*}
\int_{\mathbb{R}}
J(x+\Delta)dF_{\lfloor n\tau\rfloor +1, n}(x)
&=-\int_{\mathbb{R}}
\int_{\mathbb{R}}1_{\left\{x+\Delta\leq G(y)\right\}}H_m(y)\varphi(y)dy dF_{\lfloor n\tau\rfloor +1, n}(x)  \\
&=-\int_{\mathbb{R}}
\int_{\mathbb{R}}1_{\left\{x+\Delta\leq G(y)\right\}}dF_{\lfloor n\tau\rfloor +1, n}(x)H_m(y)\varphi(y)dy\\
&=-\int_{\mathbb{R}}
F_{\lfloor n\tau\rfloor +1, n}(G(y)-\Delta)H_m(y)\varphi(y)dy 
\end{align*}
and
\begin{align*}
\int_{\mathbb{R}}
J(x+\Delta)dF(x)
&=-\int_{\mathbb{R}}
\int_{\mathbb{R}}1_{\left\{x+\Delta\leq G(y)\right\}}H_m(y)\varphi(y)dy dF(x)  \\
&=-\int_{\mathbb{R}}
\int_{\mathbb{R}}1_{\left\{x+\Delta\leq G(y)\right\}}dF(x)H_m(y)\varphi(y)dy\\
&=-\int_{\mathbb{R}}
F(G(y)-\Delta)H_m(y)\varphi(y)dy.
\end{align*}
Regarding the difference of these terms, we obtain
\begin{align*}
&\left|\int_{\mathbb{R}}J(x+\Delta)dF_{\lfloor n\tau\rfloor +1, n}(x)-\int_{\mathbb{R}}J(x+\Delta)dF(x)\right|
\\
&=\left|\int_{\mathbb{R}}\left(F(G(y)-\Delta)-F_{\lfloor n\tau\rfloor +1, n}(G(y)-\Delta)\right)H_m(y)\varphi(y)dy\right|\\
&\leq
\int_{\mathbb{R}}\left|F(G(y)-\Delta)-F_{\lfloor n\tau\rfloor +1, n}(G(y)-\Delta)\right|\left|H_m(y)\right|\varphi(y)dy\\
&\leq
\sup\limits_{y\in \mathbb{R}}\left|F(G(y)-\Delta)-F_{\lfloor n\tau\rfloor +1, n}(G(y)-\Delta)\right|
\int_{\mathbb{R}}\left|H_m(y)\right|\varphi(y)dy,
\end{align*}
where $\int_{\mathbb{R}}\left|H_m(y)\right|\varphi(y)dy<\infty$ because of Hölder's inequality and where 
\begin{align*}
\sup_{y\in \mathbb{R}}\left|F(G(y)-\Delta)-F_{\lfloor n\tau\rfloor +1, n}(G(y)-\Delta)\right|\longrightarrow 0 \ a.s. 
\end{align*}
by \eqref{Glivenko-Cantelli}.
As a result,
$\int_{\mathbb{R}}J(x+\Delta)dF_{\lfloor n\tau\rfloor +1, n}(x)\overset{\mathcal{D}}{\longrightarrow}\int_{\mathbb{R}}J(x+\Delta)dF(x)$, so that in the end the second summand in \eqref{inequality2 Lemma1} converges to $0$ almost surely, too.

All in all, the third term on the right-hand side of 
\eqref{inequality1 Lemma1} converges to $0$ almost surely as it
is dominated by the sum of two expressions which both converge to $0$ with probability $1$.
This completes the proof of the first assertion in Lemma \ref{Lemma1}.
\hfill $\Box$

\vspace{5mm}

\begin{Kor}\label{Corollary1}
Suppose that $\left(\xi_i\right)_{i\geq 1}$  is a stationary, long-range dependent Gaussian  process   with mean $0$, variance $1$ and LRD parameter $0<D <\frac{1}{m}$, where 
 $m$ denotes the Hermite rank of the class of functions $1_{\left\{G(\xi_i)\leq x\right\}}-F(x)$, $x \in \mathbb{R}$.  
Moreover, assume that $\left(
G(\xi_i)\right)_{i\geq 1}$ has a  continuous distribution function $F$ and that $G:\mathbb{R}\longrightarrow \mathbb{R}$ is a measurable function. Then
\begin{align*}
\frac{1}{n^2}\sum\limits_{i=1}^{\lfloor n\tau\rfloor}\sum\limits_{j=\lfloor n\tau\rfloor +1}^n1_{\left\{G(\xi_i)\leq G(\xi_j)+\Delta\right\}}
\overset{P}{\longrightarrow}
\tau(1-\tau)\int_{\mathbb{R}}F(x+\Delta)dF(x)
\end{align*}
for fixed $\tau$.
\end{Kor}

\textit{Proof of Corollary \ref{Corollary1}.}
Consider the function 
$G:D\left[0, \tau\right]\longrightarrow\mathbb{R}$,
$f\mapsto f(\tau)$. As $G$
is continuous with respect to the supremum norm on $D\left[0, \tau\right]$, Corollary \ref{Corollary1} follows from Lemma \ref{Lemma1} and the continuous mapping theorem
\hfill $\Box$

\vspace{5mm}

\begin{Lem}\label{Lemma2}
Suppose that $\left(\xi_i\right)_{i\geq 1}$  is a stationary, long-range dependent Gaussian  process   with mean $0$, variance $1$ and LRD parameter $0<D <\frac{1}{m}$, where 
 $m$ denotes the Hermite rank of the class of functions $1_{\left\{G(\xi_i)\leq x\right\}}-F(x)$, $x \in \mathbb{R}$.  
Moreover, assume that $\left(
G(\xi_i)\right)_{i\geq 1}$ has a  continuous distribution function $F$ and that $G:\mathbb{R}\longrightarrow \mathbb{R}$ is a measurable function.
Then
\begin{align*}
\frac{1}{n^2}\sum\limits_{i=1}^{\lfloor n\tau\rfloor}\sum\limits_{j=\lfloor n\tau\rfloor +1}^{\lfloor n\lambda \rfloor}1_{\left\{G(\xi_i)\leq G(\xi_j)+\Delta\right\}}
\overset{P}{\longrightarrow}
\tau(\lambda-\tau)\int_{\mathbb{R}}F(x+\Delta)dF(x)
\end{align*}
for fixed $\tau$, uniformly in $\lambda\geq\tau$.
\end{Lem}

\textit{Proof of Lemma \ref{Lemma2}.}
Let $F_{k+1, t}$ denote the empirical distribution function
of $G(\xi_{k+1}), \ldots, G(\xi_t)$, i.e.
\begin{align*}
F_{k+1, t}(x)=\frac{1}{t-k}\sum\limits_{i=k+1}^t1_{\left\{G(\xi_i)\leq x\right\}}.
\end{align*}
We may therefore rewrite
\begin{align*}
\sum\limits_{i=1}^{\lfloor n\tau\rfloor}\sum\limits_{j=\lfloor n\tau\rfloor+1}^{\lfloor n\lambda \rfloor}1_{\left\{G(\xi_i)\leq G(\xi_j)+\Delta\right\}}
&=\left(\lfloor n\lambda\rfloor-\lfloor n\tau\rfloor\right)\lfloor n\tau\rfloor
\frac{1}{\left(\lfloor n\lambda\rfloor-\lfloor n\tau\rfloor\right)}\sum\limits_{j=\lfloor n\tau\rfloor +1}^{\lfloor n\lambda\rfloor}F_{\lfloor n\tau\rfloor}( G(\xi_j)+\Delta)
\\
&=\left(\lfloor n\lambda\rfloor-\lfloor n\tau\rfloor\right)\lfloor n\tau\rfloor
\int_{\mathbb{R}}F_{\lfloor n\tau\rfloor}( x+\Delta)
dF_{\lfloor n \tau \rfloor +1, \lfloor n\lambda\rfloor}(x).
\end{align*}
Furthermore,  repeated application of the triangle inequality yields
\begin{align}\label{inequality1 Lemma2}
&\sup\limits_{\tau\leq \lambda\leq 1}\left|\frac{1}{n}\left(\lfloor n\lambda\rfloor-\lfloor n\tau\rfloor\right)\int_{\mathbb{R}}F_{\lfloor n\tau\rfloor}(x+\Delta)dF_{\lfloor n\tau\rfloor +1, \lfloor n\lambda\rfloor}(x)-\left(\lambda-\tau\right)\int_{\mathbb{R}}F(x+\Delta)dF(x)\right|
\\
&\leq \sup\limits_{\tau\leq \lambda\leq 1}\left|\frac{1}{n}\left(\lfloor n\lambda\rfloor-\lfloor n\tau\rfloor\right)\int_{\mathbb{R}}\left(F_{\lfloor n\tau\rfloor}(x+\Delta)-F(x+\Delta)\right)dF_{\lfloor n\tau\rfloor +1, \lfloor n\lambda \rfloor}(x)\right|
\notag\\
&\quad \ + \sup\limits_{\tau\leq \lambda\leq 1}\left|\frac{1}{n}\left(\lfloor n\lambda\rfloor-\lfloor n\tau\rfloor\right)\int_{\mathbb{R}}F(x+\Delta)dF_{\lfloor n\tau\rfloor +1, \lfloor n\lambda \rfloor}(x)-\frac{1}{n}\left(\lfloor n\lambda\rfloor-\lfloor n\tau\rfloor\right)\int_{\mathbb{R}}F(x+\Delta)dF(x)\right|\notag\\
&\quad \ + \sup\limits_{\tau\leq \lambda\leq 1}\left|\frac{1}{n}\left(\lfloor n\lambda\rfloor-\lfloor n\tau\rfloor\right)\int_{\mathbb{R}}F(x+\Delta)dF (x)-\left(\lambda-\tau\right)\int_{\mathbb{R}}F(x+\Delta)dF (x)\right|.\notag
\end{align}

In order to prove that the stochastic process considered in Lemma \ref{Lemma2} converges to the given limit process, it is sufficient to show that the expressions on the right-hand side of the above inequality converge to $0$. We consider each of the three summands separately.

Apparently, the third term converges to $0$ since
\begin{align*}
&\sup\limits_{\tau\leq \lambda\leq 1}\left|\frac{1}{n}\left(\lfloor n\lambda\rfloor-\lfloor n\tau\rfloor\right)\int_{\mathbb{R}}F(x+\Delta)dF (x)-\left(\lambda-\tau\right)\int_{\mathbb{R}}F(x+\Delta)dF (x)\right|\\
&=\sup\limits_{\tau\leq \lambda\leq 1}\left|\frac{1}{n}\left(\lfloor n\lambda\rfloor-\lfloor n\tau\rfloor\right)-\left(\lambda-\tau\right)\right|\int_{\mathbb{R}}F(x+\Delta)dF (x)
\end{align*}
and as $\sup_{\tau\leq \lambda\leq 1}\left|\frac{1}{n}\left(\lfloor n\lambda\rfloor-\lfloor n\tau\rfloor\right)-\left(\lambda-\tau\right)\right|\longrightarrow 0$.

We have
\begin{align*}
&\sup\limits_{\tau\leq \lambda\leq 1}\left|\frac{1}{n}\left(\lfloor n\lambda\rfloor-\lfloor n\tau\rfloor\right)\int_{\mathbb{R}}\left(F_{\lfloor n\tau\rfloor}(x+\Delta)-F(x+\Delta)\right)dF_{\lfloor n\tau\rfloor +1, \lfloor n\lambda\rfloor}(x)\right|
\\
&\leq\sup\limits_{x\in \mathbb{R}}\left|F_{\lfloor n\tau\rfloor}(x+\Delta)-F(x+\Delta)\right|\sup\limits_{\tau\leq \lambda\leq 1}\left|\frac{1}{n}\left(\lfloor n\lambda\rfloor-\lfloor n\tau\rfloor\right)\right|
\end{align*}
for the first summand.
As 
$\sup_{x\in \mathbb{R}}\left|F_{\lfloor n\tau\rfloor}(x+\Delta)-F(x+\Delta)\right|$ converges to $0$ almost surely by the Glivenko-Cantelli theorem, so does the right-hand side of the above inequality.

Finally, consider the second term on the right-hand side of \eqref{inequality1 Lemma2}.
We have
\begin{align*}
&\sup\limits_{\tau\leq \lambda\leq 1}\left|\frac{1}{n}\left(\lfloor n\lambda\rfloor-\lfloor n\tau\rfloor\right)\int_{\mathbb{R}}F(x+\Delta)dF_{\lfloor n\tau\rfloor +1, \lfloor n\lambda \rfloor}(x)-\frac{1}{n}\left(\lfloor n\lambda\rfloor-\lfloor n\tau\rfloor\right)\int_{\mathbb{R}}F(x+\Delta)dF(x)\right|\\
&=\sup\limits_{\tau\leq \lambda\leq 1}\left|\frac{1}{n}\left(\lfloor n\lambda\rfloor-\lfloor n\tau\rfloor\right)\int_{\mathbb{R}}F(x+\Delta)d\left(F_{\lfloor n\tau\rfloor +1, \lfloor n\lambda \rfloor}-F\right)(x)\right|\\
&=\sup\limits_{\tau\leq \lambda\leq 1}\left|\frac{1}{n}\left(\lfloor n\lambda\rfloor-\lfloor n\tau\rfloor\right)\int_{\mathbb{R}}\left(F_{\lfloor n\tau\rfloor +1, \lfloor n\lambda\rfloor}(x)-F(x)\right)dF(x+\Delta)\right|\\
&=\frac{d_n}{n}\sup\limits_{\tau\leq \lambda\leq 1}\left|\int_{\mathbb{R}}d_n^{-1}\left(\lfloor n\lambda\rfloor-\lfloor n\tau\rfloor\right)\left(F_{\lfloor n\tau\rfloor +1, \lfloor n\lambda\rfloor}(x)-F(x)\right)-J(x)\left(Z(\lambda)-Z(\tau)\right)dF(x+\Delta)\right|\\
&\quad \ +\sup\limits_{\tau\leq \lambda\leq 1}\left|\frac{d_n}{n}\left(Z(\lambda)-Z(\tau)\right)\int_{\mathbb{R}}J(x)dF(x+\Delta)\right|\\
&\leq
\frac{d_n}{n}\sup\limits_{\tau\leq \lambda\leq 1, x\in \mathbb{R}}\left|d_n^{-1}\left(\lfloor n\lambda\rfloor-\lfloor n\tau\rfloor\right)\left(F_{\lfloor n\tau\rfloor +1, \lfloor n\lambda \rfloor}(x)-F(x)\right)-J(x)\left(Z(\lambda)-Z(\tau)\right)\right|\\
&\quad \ +\frac{d_n}{n}\sup\limits_{\tau\leq \lambda\leq 1}\left|Z(\lambda)-Z(\tau)\right|\left|\int_{\mathbb{R}}J(x)dF(x+\Delta)\right|.
\end{align*}
It follows from integration by parts that 
\begin{align*}
&\int_{\mathbb{R}}F(x+\Delta)dF_{\lfloor n\tau\rfloor +1, \lfloor n\lambda\rfloor}(x)-\int_{\mathbb{R}}F(x+\Delta)dF(x)\\
&=\int_{\mathbb{R}}F(x)dF(x+\Delta)-\int_{\mathbb{R}}F_{\lfloor n\tau\rfloor +1, \lfloor n\lambda\rfloor}(x)dF(x+\Delta).
\end{align*}
Furthermore,
\begin{align*}
\left(\lfloor n\lambda\rfloor-\lfloor n\tau\rfloor\right)\left(F_{\lfloor n\tau\rfloor +1, \lfloor n\lambda \rfloor}(x)-F(x)\right)
&=\sum\limits_{i=\lfloor n\tau\rfloor +1}^{\lfloor n\lambda\rfloor}1_{\{G(\xi_i)\leq x\}}-\left(\lfloor n\lambda\rfloor-\lfloor n\tau\rfloor\right)F(x)
\\
&=\lfloor n\lambda\rfloor F_{\lfloor n\lambda\rfloor}(x)-\lfloor n\tau\rfloor F_{\lfloor n\tau\rfloor}(x)-\lfloor n\lambda\rfloor F(x)+\lfloor n\tau\rfloor F(x)\\
&=\lfloor n\lambda\rfloor \left(F_{\lfloor n\lambda\rfloor}(x)-F(x)\right)-\lfloor n\tau\rfloor \left(F_{\lfloor n\tau\rfloor}(x)- F(x)\right).
\end{align*}
As a result, 
\begin{align*}
&\sup\limits_{\tau\leq \lambda\leq 1, x\in \mathbb{R}}\left|d_n^{-1}\left(\lfloor n\lambda\rfloor-\lfloor n\tau\rfloor\right)\left(F_{\lfloor n\tau\rfloor +1, \lfloor n\lambda \rfloor}(x)-F(x)\right)-J(x)\left(Z(\lambda)-Z(\tau)\right)\right|\\
&=\sup\limits_{\tau\leq \lambda\leq 1, x\in \mathbb{R}}\left|d_n^{-1}\left(\lfloor n\lambda\rfloor \left(F_{\lfloor n\lambda\rfloor}(x)-F(x)\right)-\lfloor n\tau\rfloor \left(F_{\lfloor n\tau\rfloor}(x)- F(x)\right)\right)-J(x)\left(Z(\lambda)-Z(\tau)\right)\right|\\
&\leq \sup\limits_{\tau\leq \lambda\leq 1, x\in \mathbb{R}}\left|d_n^{-1}\lfloor n\lambda\rfloor \left(F_{\lfloor n\lambda\rfloor}(x)-F(x)\right)-J(x)Z(\lambda)\right|\\
&\quad \ +\sup\limits_{\tau\leq \lambda\leq 1, x\in \mathbb{R}}\left|d_n^{-1}\lfloor n\tau\rfloor \left(F_{\lfloor n\tau\rfloor}(x)- F(x)\right)-J(x)Z(\tau)\right|.
\end{align*}
Again, we may assume 
without loss of generality   that
\begin{align*}
\sup\limits_{\lambda, x}\left|d_n^{-1}\lfloor n\lambda\rfloor\left(F_{\lfloor n\lambda\rfloor}(x)-F(x)\right)-J(x)Z(\lambda)\right|\longrightarrow 0
\end{align*}
almost surely,  as pointed out in the proof of Lemma \ref{Lemma1}. 
Since $\frac{d_n}{n}\longrightarrow 0$  by definition of $d_n$, we may conclude that the third summand on the right hand side of \eqref{inequality1 Lemma2}  converges to $0$, too. This completes the proof of Lemma \ref{Lemma2}.
\hfill $\Box$

\vspace{5mm}

\textit{Proof of Theorem \ref{consistency}.}
We have
\begin{align*}
T_n(\tau_1, \tau_2)&=\sup\limits_{k\in \left[\lfloor n\tau_1\rfloor, \lfloor n\tau_2\rfloor \right]}G_n(k)\\
&\geq G_n(k^*),
\end{align*}
where 
\begin{align*}
G_n(k^*)=\frac{\left|\sum\limits_{i=1}^{k^*}\sum\limits_{j=k^*+1}^n\left(1_{\left\{X_i\leq X_j\right\}}-\frac{1}{2}\right)\right|}{\left\{\frac{1}{n}\sum\limits_{t=1}^{k^*}S_t^2(1, k^*)+\frac{1}{n}\sum\limits_{t=k^*+1}^nS_t^2(k^*+1, n)\right\}^{\frac{1}{2}}}
\end{align*}
and where $k^*=\lfloor n\tau\rfloor$ denotes the location of the change-point.
Thus, 
it suffices to show that  $G_n(k^*)\overset{P}{\longrightarrow}\infty$.
For this purpose, we rewrite
\begin{align*}
G_n(k^*)=\frac{\frac{1}{n^2}\left|\sum\limits_{i=1}^{k^*}\sum\limits_{j=k^*+1}^n\left(1_{\left\{X_i\leq X_j\right\}}-\frac{1}{2}\right)\right|}{\frac{1}{n^2}\left\{\frac{1}{n}\sum\limits_{t=1}^{k^*}S_t^2(1, k^*)+\frac{1}{n}\sum\limits_{t=k^*+1}^nS_t^2(k^*+1, n)\right\}^{\frac{1}{2}}}.
\end{align*}

We will prove that the numerator of $G_n(k^*)$ converges to a positive constant, whereas the denominator tends to $0$ in order to show divergence to $\infty$.

First, we turn to the denominator, which equals
\begin{align*}
\frac{1}{n^2}\left\{\int_0^{\tau}S_{\lfloor nr\rfloor}^2(1, k^*)dr+\int_{\tau}^1S_{\lfloor nr\rfloor}^2(k^*+1, n)dr\right\}^{\frac{1}{2}}.
\end{align*}
Note that for $i\leq k^*$
\begin{align*}
\sum\limits_{j=1}^n1_{\left\{X_i\leq X_j\right\}}
&=\sum\limits_{j=1}^{k^*}1_{\left\{\mu+G(\xi_i)\leq \mu+G(\xi_j)\right\}}+\sum\limits_{j=k^*+1}^n1_{\left\{\mu+G(\xi_i)\leq \mu+G(\xi_j)+\Delta\right\}}\\
&=\sum\limits_{j=1}^n1_{\left\{G(\xi_i)\leq G(\xi_j)\right\}}+\sum\limits_{j=k^*+1}^n1_{\left\{G(\xi_j)<G(\xi_i)\leq G(\xi_j)+\Delta\right\}}.
\end{align*}
Therefore,
\begin{align*}
S_t(1, k^*)
&=-\sum\limits_{h=1}^t\left(
\sum\limits_{j=1}^n1_{\left\{X_h\leq X_j\right\}}-\frac{1}{k^*}\sum\limits_{i=1}^{k^*}\sum\limits_{j=1}^n1_{\left\{X_i\leq X_j\right\}}
\right)\\
&=-\sum\limits_{h=1}^t\left(
\sum\limits_{j=1}^n1_{\left\{G(\xi_h)\leq G(\xi_j)\right\}}-\frac{1}{k^*}\sum\limits_{i=1}^{k^*}\sum\limits_{j=1}^n1_{\left\{G(\xi_i)\leq G(\xi_j)\right\}}
\right)\\
&\quad \ -\sum\limits_{h=1}^t\left(
\sum\limits_{j=k^*+1}^n1_{\left\{G(\xi_j)<G(\xi_h)\leq G(\xi_j)+\Delta\right\}}-\frac{1}{k^*}\sum\limits_{i=1}^{k^*}\sum\limits_{j=k^*+1}^n1_{\left\{G(\xi_j)<G(\xi_i)\leq G(\xi_j)+\Delta\right\}}
\right).
\end{align*}
We treat the expression 
 $S_t(1, k^*)$ as 
 sum of the following terms
\begin{align*}
&\hat{S}_t(1, k^*)=-\sum\limits_{h=1}^t\left(
\sum\limits_{j=1}^n1_{\left\{G(\xi_h)\leq G(\xi_j)\right\}}-\frac{1}{k^*}\sum\limits_{i=1}^{k^*}\sum\limits_{j=1}^n1_{\left\{G(\xi_i)\leq G(\xi_j)\right\}}
\right),\\
&\tilde{S}_t(1, k^*)=-\sum\limits_{h=1}^t\left(
\sum\limits_{j=k^*+1}^n1_{\left\{G(\xi_j)<G(\xi_h)\leq G(\xi_j)+\Delta\right\}}-\frac{1}{k^*}\sum\limits_{i=1}^{k^*}\sum\limits_{j=k^*+1}^n1_{\left\{G(\xi_j)<G(\xi_i)\leq G(\xi_j)+\Delta\right\}}
\right).
\end{align*}

For the first summand we get
\begin{align*}
\hat{S}_t(1, k^*)
&=-\sum\limits_{h=1}^t\left(
\sum\limits_{j=1}^n1_{\left\{G(\xi_h)\leq G(\xi_j)\right\}}-\frac{1}{k^*}\sum\limits_{i=1}^{k^*}\sum\limits_{j=1}^n1_{\left\{G(\xi_i)\leq G(\xi_j)\right\}}
\right)\\
&=-\sum\limits_{i=1}^t
\sum\limits_{j=1}^t1_{\left\{G(\xi_i)\leq G(\xi_j)\right\}}-\sum\limits_{i=1}^t\sum\limits_{j=t+1}^n1_{\left\{G(\xi_i)\leq G(\xi_j)\right\}}\\
&\quad \ +
\frac{t}{k^*}
\sum\limits_{i=1}^{k^*}\sum\limits_{j=1}^{k^*}1_{\left\{G(\xi_i)\leq G(\xi_j)\right\}} +\frac{t}{k^*}\sum\limits_{i=1}^{k^*}\sum\limits_{j=k^*+1}^n1_{\left\{G(\xi_i)\leq G(\xi_j)\right\}}\\
&=-\frac{t(t+1)}{2}-\sum\limits_{i=1}^t\sum\limits_{j=t+1}^n1_{\left\{G(\xi_i)\leq G(\xi_j)\right\}}+
\frac{t}{k^*}
\frac{k^*(k^*+1)}{2} +\frac{t}{k^*}\sum\limits_{i=1}^{k^*}\sum\limits_{j=k^*+1}^n1_{\left\{G(\xi_i)\leq G(\xi_j)\right\}}\\
&=-\frac{t^2}{2}-\sum\limits_{i=1}^t\sum\limits_{j=t+1}^n1_{\left\{G(\xi_i)\leq G(\xi_j)\right\}}+
\frac{tk^*}{2} +\frac{t}{k^*}\sum\limits_{i=1}^{k^*}\sum\limits_{j=k^*+1}^n1_{\left\{G(\xi_i)\leq G(\xi_j)\right\}}.
\end{align*}
We have
\begin{align*}
&\frac{1}{n^2}\sum\limits_{i=1}^{\lfloor n\lambda\rfloor}\sum\limits_{j=\lfloor n\lambda\rfloor+1}^n1_{\left\{G(\xi_i)\leq G(\xi_j)\right\}}
\overset{P}{\longrightarrow}
\frac{\lambda(1-\lambda)}{2}
\end{align*}
uniformly in $\lambda$ because
\begin{align*}
\frac{1}{nd_n}\sum\limits_{i=1}^{\lfloor n\lambda\rfloor}\sum\limits_{j=\lfloor n\lambda\rfloor+1}^n\left(1_{\left\{G(\xi_i)\leq G(\xi_j)\right\}}-\frac{1}{2}\right)
\overset{\mathcal{D}}{\longrightarrow}
\left(Z(\lambda)-\lambda Z(1)\right)\int_{\mathbb{R}}J(x)dF(x)
\end{align*}
uniformly in $\lambda$ by Theorem 1.1 in \cite{DehlingRoochTaqqu2013a} and as $\frac{d_n}{n}\longrightarrow 0$. 
We may conclude from this and Corollary \ref{Corollary1} that
$\frac{1}{n^2}\hat{S}_{\lfloor n\lambda\rfloor}(1, \lfloor k^*\rfloor)
\overset{P}{\longrightarrow}0$
uniformly in $\lambda\leq \tau$.

Because of
\begin{align*}
1_{\left\{G(\xi_j)<G(\xi_i)\leq G(\xi_j)+\Delta\right\}}
&=
1_{\left\{G(\xi_i)\leq G(\xi_j)+\Delta\right\}}-1_{\left\{G(\xi_i)\leq G(\xi_j)\right\}},
\end{align*}
the second summand  can be written as
\begin{align*}
\tilde{S}_{t}(1, {k^*})
&=-\sum\limits_{h=1}^{t}\left(
\sum\limits_{j={k^*}+1}^n1_{\left\{G(\xi_j)<G(\xi_h)\leq G(\xi_j)+\Delta\right\}}-\frac{1}{{k^*}}\sum\limits_{i=1}^{{k^*}}\sum\limits_{j={k^*}+1}^n1_{\left\{G(\xi_j)<G(\xi_i)\leq G(\xi_j)+\Delta\right\}}
\right)\\
&=-\sum\limits_{i=1}^{t}
\sum\limits_{j={k^*}+1}^n1_{\left\{G(\xi_i)\leq G(\xi_j)+\Delta\right\}}+\frac{t}{{k^*}}\sum\limits_{i=1}^{{k^*}}\sum\limits_{j={k^*}+1}^n1_{\left\{G(\xi_i)\leq G(\xi_j)+\Delta\right\}}
\\
&\quad \ +\sum\limits_{i=1}^{t}
\sum\limits_{j={k^*}+1}^n1_{\left\{G(\xi_i)\leq G(\xi_j)\right\}}-\frac{t}{{k^*}}\sum\limits_{i=1}^{{k^*}}\sum\limits_{j={k^*}+1}^n1_{\left\{G(\xi_i)\leq G(\xi_j)\right\}}.
\end{align*}
Due to Lemma \ref{Lemma1} and  Corollary \ref{Corollary1},
$\frac{1}{n^2}\tilde{S}_{\lfloor n\lambda\rfloor}(1, k^*)$
converges in probability to $0$, as well.

All in all, the previous considerations yield 
\begin{align*}
\int_{0}^{\tau}\left\{\frac{1}{n^2}
S_{\lfloor nr\rfloor}(1, k^*)\right\}^2dr\overset{P}{\longrightarrow} 0
\end{align*}
as $G:D\left[0, \tau\right]\longrightarrow \mathbb{R}$, $f\mapsto \int_{0}^{\tau}\left(f(s)\right)^2ds$, is   continuous with respect to the  supremum norm on $D\left[0, \tau\right]$.

 In analogy to the previous argumentation it can be shown that 
\begin{align*}
\int_{\tau}^1\left\{\frac{1}{n^2}S_{\lfloor nr\rfloor}(k^*+1, n)\right\}^2dr \overset{P}{\longrightarrow}0.
\end{align*}

For this purpose, note that, if $i> {k^*}$,
\begin{align*}
\sum\limits_{j=1}^n1_{\left\{X_i\leq X_j\right\}}
&=\sum\limits_{j=1}^{k^*}1_{\left\{\mu+G(\xi_i)+\Delta\leq \mu+G(\xi_j)\right\}}+\sum\limits_{j={k^*}+1}^n1_{\left\{\mu+G(\xi_i)+\Delta\leq \mu+G(\xi_j)+\Delta\right\}}\\
&=\sum\limits_{j=1}^n1_{\left\{G(\xi_i)\leq G(\xi_j)\right\}}-\sum\limits_{j=1}^{k^*}1_{\{G(\xi_i)\leq G(\xi_j)<G(\xi_i)+\Delta\}}.
\end{align*}
Therefore,
\begin{align*}
S_t({k^*}+1, n)
&=-\sum\limits_{h={k^*}+1}^t\left(
\sum\limits_{j=1}^n1_{\left\{X_h\leq X_j\right\}}-\frac{1}{n-{k^*}}\sum\limits_{i={k^*}+1}^n\sum\limits_{j=1}^n1_{\left\{X_i\leq X_j\right\}}
\right)\\
&=-\sum\limits_{h={k^*}+1}^t\left(
\sum\limits_{j=1}^n1_{\left\{G(\xi_h)\leq G(\xi_j)\right\}}-\frac{1}{n-{k^*}}\sum\limits_{i={k^*}+1}^n\sum\limits_{j=1}^n1_{\left\{G(\xi_i)\leq G(\xi_j)\right\}}
\right)\\
&\quad \ +\sum\limits_{h={k^*}+1}^t\left(
\sum\limits_{j=1}^{k^*}1_{\{G(\xi_h)\leq G(\xi_j)<G(\xi_h)+\Delta\}}-\frac{1}{n-{k^*}}\sum\limits_{i={k^*}+1}^n\sum\limits_{j=1}^{k^*}1_{\{G(\xi_i)\leq G(\xi_j)<G(\xi_i)+\Delta\}}
\right).
\end{align*}
Hence, we consider 
 $S_t({k^*}+1, n)$ as sum of the  expressions below
\begin{align*}
&\hat{S}_t({k^*}+1, n)=-\sum\limits_{h={k^*}+1}^t\left(
\sum\limits_{j=1}^n1_{\left\{G(\xi_h)\leq G(\xi_j)\right\}}-\frac{1}{n-{k^*}}\sum\limits_{i={k^*}+1}^n\sum\limits_{j=1}^n1_{\left\{G(\xi_i)\leq G(\xi_j)\right\}}
\right),\\
&\tilde{S}_t({k^*}+1, n)=\sum\limits_{h={k^*}+1}^t\left(
\sum\limits_{j=1}^{k^*}1_{\left\{G(\xi_h)\leq G(\xi_j)< G(\xi_h)+\Delta\right\}}-\frac{1}{n-{k^*}}\sum\limits_{i={k^*}+1}^n\sum\limits_{j=1}^{k^*}1_{\left\{G(\xi_i)\leq G(\xi_j)< G(\xi_i)+\Delta\right\}}
\right).
\end{align*}
The following representation arises from rather simple transformations
\begin{align*}
\hat{S}_t({k^*}+1, n)&=-\sum\limits_{h={k^*}+1}^t\left(\sum\limits_{j=1}^n1_{\{G(\xi_h)\leq G(\xi_j)\}}-\frac{1}{n-{k^*}}\sum\limits_{i={k^*}+1}^n\sum\limits_{j=1}^n1_{\{G(\xi_i)\leq G(\xi_j)\}}\right)\\
&=-\sum\limits_{i=1}^t\sum\limits_{j=t+1}^n1_{\{G(\xi_i)\leq G(\xi_j)\}}-\sum\limits_{i=1}^t\sum\limits_{j=1}^t1_{\{G(\xi_i)\leq G(\xi_j)\}}\\
&\quad \ +\sum\limits_{i=1}^{k^*}\sum\limits_{j={k^*}+1}^n1_{\{G(\xi_i)\leq G(\xi_j)\}}+
\sum\limits_{i=1}^{k^*}\sum\limits_{j=1}^{k^*}1_{\{G(\xi_i)\leq G(\xi_j)\}}\\
&\quad \ +\frac{t-{k^*}}{n-{k^*}}\sum\limits_{i={k^*}+1}^n\sum\limits_{j=1}^{k^*}1_{\{G(\xi_i)\leq G(\xi_j)\}}+
\frac{t-{k^*}}{n-{k^*}}\sum\limits_{i={k^*}+1}^n\sum\limits_{j={k^*}+1}^n1_{\{G(\xi_i)\leq G(\xi_j)\}}\\
&=-\sum\limits_{i=1}^t\sum\limits_{j=t+1}^n1_{\{G(\xi_i)\leq G(\xi_j)\}}-\frac{t(t+1)}{2} +\sum\limits_{i=1}^{k^*}\sum\limits_{j={k^*}+1}^n1_{\{G(\xi_i)\leq G(\xi_j)\}}+
\frac{{k^*}({k^*}+1)}{2}\\
&\quad \ +\frac{t-{k^*}}{n-{k^*}}\sum\limits_{i={k^*}+1}^n\sum\limits_{j=1}^{k^*}\left(1-1_{\{G(\xi_j)\leq G(\xi_i)\}}\right)+
\frac{t-{k^*}}{n-{k^*}}\frac{(n-{k^*})(n-{k^*}+1)}{2}\\
&=-\sum\limits_{i=1}^t\sum\limits_{j=t+1}^n1_{\{G(\xi_i)\leq G(\xi_j)\}}-\frac{t(t+1)}{2} +\sum\limits_{i=1}^{k^*}\sum\limits_{j={k^*}+1}^n1_{\{G(\xi_i)\leq G(\xi_j)\}}+
\frac{{k^*}({k^*}+1)}{2}\\
&\quad \ +(t-{k^*}){k^*}-\frac{t-{k^*}}{n-{k^*}}\sum\limits_{j=1}^{k^*}\sum\limits_{i={k^*}+1}^n1_{\{G(\xi_j)\leq G(\xi_i)\}}+
\frac{(t-{k^*})(n-{k^*}+1)}{2}.
\end{align*}
Based on Lemma \ref{Lemma1} and Corollary \ref{Corollary1}, the  argumentation that
also established 
 $\frac{1}{n^2}\hat{S}_{\lfloor n\lambda\rfloor}(1, k^*)\overset{P}{\longrightarrow}0$
 yields
$\frac{1}{n^2}\hat{S}_{\lfloor n\lambda\rfloor}({k^*}+1, n)
\overset{P}{\longrightarrow}0$ uniformly in $\lambda\geq \tau$.

Likewise, it can be shown that 
$\frac{1}{n^2}\tilde{S}_{\lfloor n\lambda\rfloor}({k^*}+1, n)
\overset{P}{\longrightarrow}0$.
First of all, we note that
\begin{align*}
1_{\left\{G(\xi_i)\leq G(\xi_j)< G(\xi_i)+\Delta\right\}}
&=
1_{\left\{G(\xi_j)\leq G(\xi_i)+\Delta\right\}}-1_{\left\{G(\xi_j)\leq G(\xi_i)\right\}}
\end{align*}
almost surely if $i\neq j$.
Thereby,
\begin{align*}
\tilde{S}_t({k^*}+1, n)
&=\sum\limits_{h={k^*}+1}^t\left(
\sum\limits_{j=1}^{k^*}1_{\{G(\xi_h)\leq G(\xi_j)< G(\xi_h)+\Delta\}}-\frac{1}{n-{k^*}}\sum\limits_{i={k^*}+1}^n\sum\limits_{j=1}^{k^*}1_{\{G(\xi_i)\leq G(\xi_j)< G(\xi_i)+\Delta\}}
\right)\\
&=
\sum\limits_{j=1}^{k^*} \sum\limits_{i={k^*}+1}^t1_{\left\{G(\xi_j)\leq G(\xi_i)+\Delta\right\}}-\frac{t-{k^*}}{n-{k^*}}\sum\limits_{j=1}^{k^*}\sum\limits_{i={k^*}+1}^n1_{\left\{G(\xi_j)\leq G(\xi_i)+\Delta\right\}}
\\
&\quad \ -
\sum\limits_{j=1}^{k^*}\sum\limits_{i={k^*}+1}^t1_{\left\{G(\xi_j)\leq G(\xi_i)\right\}}+\frac{t-{k^*}}{n-{k^*}}\sum\limits_{j=1}^{k^*}\sum\limits_{i={k^*}+1}^n1_{\left\{G(\xi_j)\leq G(\xi_i)\right\}}.
\end{align*}
As a result, we have $\frac{1}{n^2}\tilde{S}_{\lfloor n\lambda\rfloor}({k^*}+1, n)
\overset{P}{\longrightarrow}0$ by Lemma \ref{Lemma2} and Corollary \ref{Corollary1}.

As both terms, $\frac{1}{n^2}\hat{S}_{\lfloor n\lambda\rfloor}({k^*}+1, n)$ as well as $\frac{1}{n^2}\tilde{S}_{\lfloor n\lambda\rfloor}({k^*}+1, n)$, converge in probability to $0$ uniformly in $\lambda\geq \tau$, 
it follows that
\begin{align*}
\int_{\tau}^{1}\left\{\frac{1}{n^2}
S_{\lfloor n r\rfloor}(k^*+1, n)\right\}^2dr \overset{P}{\longrightarrow}0.
\end{align*}

On the basis of the previous considerations 
we may conclude that the denominator of $G_n(k^*)$
converges in probability to $0$.

In order to prove the consistency of the self-normalized Wilcoxon change-point test,
it therefore remains to show that the numerator of $G_n(k^*)$, given by
\begin{align*}
\frac{1}{n^2}\left|\sum\limits_{i=1}^{k^*}\sum\limits_{j=k^*+1}^n\left(1_{\left\{X_i\leq X_j\right\}}-\frac{1}{2}\right)\right|,
\end{align*} converges to 
a non-negative constant.

We have
\begin{align*}
\frac{1}{n^2}\left|\sum\limits_{i=1}^{k^*}\sum\limits_{j=k^*+1}^n\left(1_{\left\{X_i\leq X_j\right\}}-\frac{1}{2}\right)\right|
&=\left|\frac{1}{n^2}\sum\limits_{i=1}^{k^*}\sum\limits_{j=k^*+1}^n1_{\left\{G(\xi_i)\leq G(\xi_j)+\Delta\right\}}-\frac{1}{n^2}\frac{k^*(n-k^*)}{2}\right|.
\end{align*}
Therefore,
\begin{align}\label{numerator}
\frac{1}{n^2}\left|\sum\limits_{i=1}^{k^*}\sum\limits_{j=k^*+1}^n\left(1_{\left\{X_i\leq X_j\right\}}-\frac{1}{2}\right)\right|\overset{P}{\longrightarrow}
\tau(1-\tau)\int_{\mathbb{R}}\left(F(x+\Delta)-F(x)\right)dF(x)
\end{align}
by  Corollary \ref{Corollary1} and since
$\frac{1}{n^2}\frac{k^*(n-k^*)}{2}\longrightarrow \frac{\tau(1-\tau)}{2}$.
As the limit in \eqref{numerator} does not vanish,  $G_n(k^*)$ diverges to $\infty$ and we
thus have proved Theorem
\ref{consistency}.
\hfill $\Box$

\vspace{5mm}

\textit{Proof of Theorem \ref{asymptotic distribution under local A}.}
Note that because of the corresponding sample path properties of the stochastic process $Z$, the sample paths 
of
\begin{align*}
W_{m, \tau}^{*}(\lambda)=\left(Z(\lambda)-\lambda Z(1)\right)\int_{\mathbb{R}}J(x)dF(x)+c\delta_{\tau}(\lambda)\int_{\mathbb{R}}f^2(x)dx, \  0\leq \lambda\leq 1,
\end{align*}
are almost surely continuous and nowhere differentiable. 

The same argument as in the proof of Theorem \ref{asymptotic  distribution under H}
shows that $T_n(\tau_1, \tau_2)$ converges in distribution to  
\begin{align*}
\sup\limits_{\lambda\in \left[\tau_1, \tau_2\right]}\frac{\left|W_{m, \tau}^{*}(\lambda)\right|}{\left\{\int_0^{\lambda}\left(W_{m, \tau}^{*}(r)-\frac{r}{\lambda}W_{m, \tau}^{*}(\lambda)\right)^2dr+\int_{\lambda}^1\left(W_{m, \tau}^{*}(r)-\frac{1-r}{1-\lambda}W_{m, \tau}^{*}(\lambda)\right)^2dr\right\}^{\frac{1}{2}}}.
\end{align*}

The numerator of the limit process equals
\begin{align*}
\left|\int_{\mathbb{R}}J(x)dF(x)\left(Z(\lambda)-\lambda Z(1)\right)+c\delta_{\tau}(\lambda)\int_{\mathbb{R}}f^2(x)dx\right|.
\end{align*}
Moreover, for the quantities in the denominator it holds that
\begin{align*}
W_{m, \tau}^*(r)-\frac{r}{\lambda}W_{m, \tau}^*(\lambda)
&=
\left(Z(r)-r Z(1)\right)\int_{\mathbb{R}}J(x)dF(x)+c\delta_{\tau}
(r)\int_{\mathbb{R}}f^2(x)dx\\
&\quad \ -\frac{r}{\lambda}\left(\left(Z(\lambda)-\lambda Z(1)\right)\int_{\mathbb{R}}J(x)dF(x)+c\delta_{\tau}(\lambda)\int_{\mathbb{R}}f^2(x)dx\right)\\
&=\int_{\mathbb{R}}J(x)dF(x)
\left(Z(\lambda)-\frac{r}{\lambda}Z(\lambda)\right)+c\int_{\mathbb{R}}f^2(x)dx\left(\delta_{\tau}(r)-\frac{r}{\lambda}\delta_{\tau}(\lambda)\right)
\end{align*}
and
\begin{align*}
W_{m, \tau}^*(r)-\frac{1-r}{1-\lambda}W_{m, \tau}^*(\lambda)
&=\left(Z(r)-r Z(1)\right)\int_{\mathbb{R}}J(x)dF(x)+c\delta_{\tau}
(r)\int_{\mathbb{R}}f^2(x)dx\\
&\quad \ -\frac{1-r}{1-\lambda}\left\{\left(Z(\lambda)-\lambda Z(1)\right)\int_{\mathbb{R}}J(x)dF(x)+c\delta_{\tau}
(\lambda)\int_{\mathbb{R}}f^2(x)dx\right\}\\
&=\int_{\mathbb{R}}J(x)dF(x)\left\{Z(r)-rZ(1)-\frac{1-r}{1-\lambda}\left(Z(\lambda)-\lambda Z(1)\right)\right\}\\
&\quad \  +c\int_{\mathbb{R}}f^2(x)dx
\left(\delta_{\tau}(r)-\frac{1-r}{1-\lambda}\delta_{\tau}(\lambda)\right)
.
\end{align*}
\hfill $\Box$

\vspace{5mm}

\thispagestyle{plain}  
\bibliography{Paper}

\end{document}